\documentclass[11pt]{article}
\usepackage{amsmath,amssymb}
\usepackage{amsthm}
\usepackage{amsfonts}
\usepackage{natbib}
\usepackage{graphicx}
\usepackage{relsize}
\usepackage{catchfilebetweentags}
\usepackage{layout}
\usepackage{bm}
\usepackage{caption}
\usepackage{float}

\usepackage{lscape} 
\usepackage{pdflscape}

\usepackage{multirow}
\usepackage{multicol}
\usepackage{bigints}
\usepackage[margin=1in]{geometry}
\usepackage[onehalfspacing]{setspace}
\usepackage{mathtools}

\usepackage{rotating}

\usepackage[shortlabels]{enumitem} 

\allowdisplaybreaks

\newtheorem{theorem}{Theorem}
\newtheorem{corollary}{Corollary}
\newtheorem{remark}{Remark}
\newtheorem{assumption}{Assumption}
\newtheorem{lemma}{Lemma}

\newcommand{\m}[1]{\mathcal{#1}}

\newcommand{\boup}[1]{\boldsymbol{\mathrm{#1}}}

\newcommand{\ha}[1]{\widehat{#1}}
\newcommand{\ti}[1]{\widetilde{#1}}

\newcommand{\mme}{\mathbb{E}}
\newcommand{\mmi}{\mathbb{I}}
\newcommand{\mmm}{\mathbb{M}}
\newcommand{\mmn}{\mathbb{N}}
\newcommand{\mmp}{\mathbb{P}}
\newcommand{\mmr}{\mathbb{R}}

\newcommand{\mmw}{\mathbb{W}}
\newcommand{\mmx}{\mathbb{X}}

\newcommand{\eps}{\varepsilon}

\newcommand{\dto}{\overset{d}{\to}} 
\newcommand{\pto}{\overset{p}{\to}} 

\numberwithin{equation}{section}

\graphicspath{{figures_tables/}} 

\makeatletter
\def\input@path{{figures_tables/}}
\makeatother

\begin{document}

\title{Impossible Inference in Econometrics: Theory and Applications}

\author{
\renewcommand{\thefootnote}{\alph{footnote}}
Marinho Bertanha\footnotemark[1]
\and 
\renewcommand{\thefootnote}{\alph{footnote}}
Marcelo J. Moreira\footnotemark[2]
}

\date{\vspace{1cm} This version: September 24, 2019 
\\
First version: October 11, 2016}

{
\renewcommand{\thefootnote}{\alph{footnote}}
\footnotetext[1]{
Department of Economics, University of Notre Dame. Address: 3060 Jenkins Nanovic Halls, Notre Dame, IN 46556, USA.
Email: mbertanha@nd.edu. Website: www.nd.edu/$\sim$mbertanh.
  }
\footnotetext[2]{Corresponding author.
FGV EPGE. Address: Praia de Botafogo, 190 11th floor, Rio de Janeiro - RJ 22250-040, Brazil.
Email: mjmoreira@fgv.br. Website: epge.fgv.br/en/professor/marcelo-moreira.
}  
}

\maketitle

{
\begin{abstract}
This paper studies models in which hypothesis tests have trivial power, that is, power smaller than size.
This testing impossibility, or impossibility type A, arises when any alternative is \textit{not distinguishable} from the null.
We also study settings in which it is impossible to have almost surely bounded confidence sets for a parameter of interest.
This second type of impossibility (type B) occurs under a condition weaker than the condition for type A impossibility: the parameter of interest must be \textit{nearly unidentified}.
Our theoretical framework connects many existing publications on impossible inference that rely on different notions of topologies to show models are not distinguishable or nearly unidentified.
We also derive both types of impossibility using the weak topology induced by convergence in distribution. 
Impossibility in the weak topology is often easier to prove, it is applicable for many widely-used tests, and it is
useful for robust hypothesis testing.
We conclude by demonstrating impossible inference in multiple economic applications of models with discontinuity and time-series models. 

\end{abstract}
}

\vspace{1cm}

{



}

\newpage

\section{Introduction}

\indent

The goal of most empirical studies is to estimate parameters of a population
statistical model using a random sample of data. The difference between
estimates and population parameters is uncertain because sample data do not
have all information about the population. Statistical inference
provides methods for quantifying this uncertainty. Typical approaches
include hypothesis testing and confidence sets. In a hypothesis test, the
researcher divides all possible population models into two sets of models.
The null set includes models which the researcher suspects to be
false. The alternative set includes all other likely models. It is desirable
to control the size of the test, that is, the error probability of rejecting
the null set when the null set contains the true model. A powerful test has
a small error probability of failing to reject the null set when the true model
is outside the null set. Another approach is to use the data to build a
confidence set for the unknown value of parameters of the true model. 
The
researcher needs to control the error probability that the confidence set
excludes the true value. 
Error probabilities must be controlled uniformly over the entire set of likely models.
This paper studies necessary and sufficient conditions for the impossibility of controlling error probabilities of hypothesis tests and confidence sets.

Previous work demonstrates the impossibility of controlling error probabilities of tests and confidence sets in specific settings. There are essentially two types of impossibility found in the literature.
The first type of impossibility says that any hypothesis test has power limited by size.
That is, it is impossible to find a powerful test that controls size.
We call this impossibility type A.
The second type of impossibility states that any confidence set that is almost surely (a.s.) bounded 
has error probability arbitrarily equal to one (i.e. zero confidence level).  
In other words, it is impossible for finite bounds to contain the true value of parameters with high probability.
We call this impossibility type B.
Despite being related, both types of impossibility often appear disconnected in the existing literature.

\begin{sidewaysfigure}
\caption{Impossibility Literature Diagram}
\label{diagram}
\begin{center}
\includegraphics[trim={.75in .45in 0 0},clip,width=9.5in]{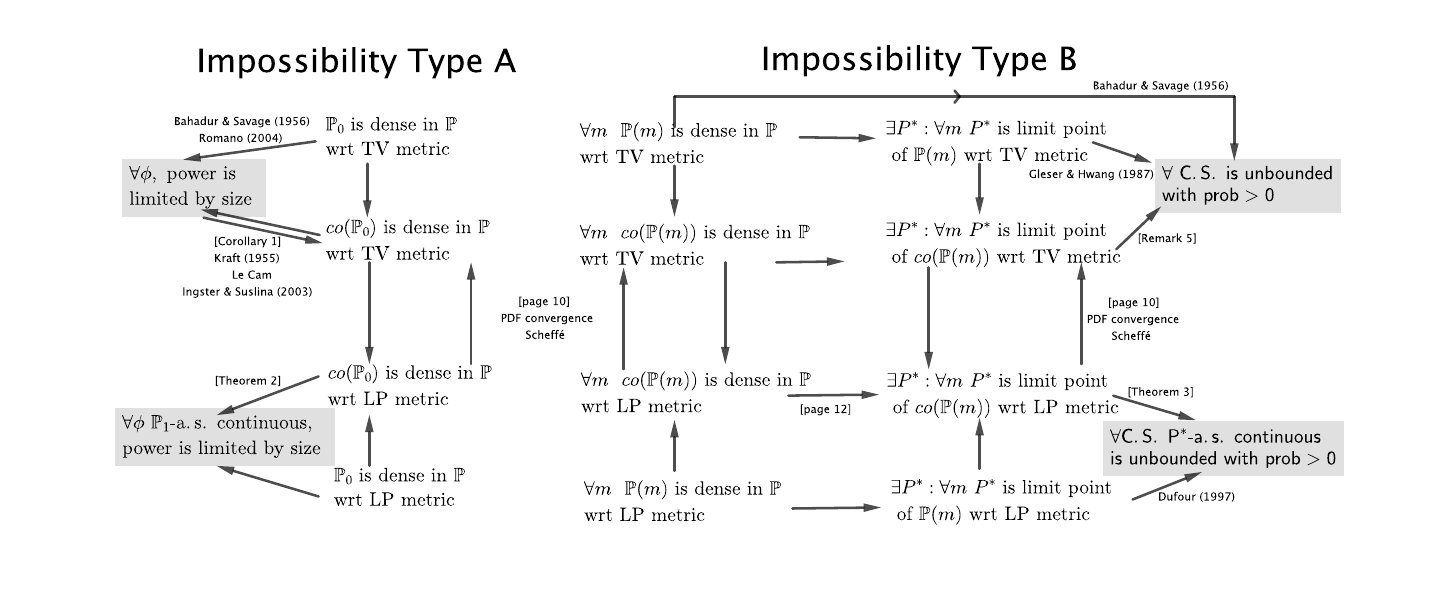}
\end{center}
\par
\caption*{\footnotesize Notes: 
the diagram illustrates the relationships between the different versions of impossibility found in the literature. 
Arrows without labels or arrows with references in square brackets are relationships made explicit by this paper.
Impossibility type A says that every test function $\phi$ has maximum power less than or equal to size.
The set of all likely models is $\mmp$, which is the union of models under the null $\mmp_0$
 and alternative hypothesis $\mmp_1$.
 A test $\phi$ is said to be $\mmp_1$-a.s. continuous if the set of discontinuity points of $\phi$
 has zero probability under every $Q \in \mmp_1$.
The set $co(\mmp_0)$ denotes the convex hull of $\mmp_0$.
The TV metric is the Total Variation metric (Equation (\ref{def:tv-metric}) in Section \ref{sec:theo}). 
The LP metric is the L\'{e}vy-Prokhorov metric (Equation (\ref{def:lp-metric}) in Section \ref{sec:theo}). 
The set $\mmp_0$ is dense in $\mmp$ wrt a metric $d(\cdot,\cdot)$
if, for every $Q \in \mmp_1$, there exists a sequence $\{ P_k\}_k \subseteq \mmp_0$
such that $d(P_k,Q) \to 0$.
Impossibility type B says that every confidence set (C.S.) is unbounded with positive probability for some distributions in $\mmp$. 
The subset $\mmp(m)$ denotes all models $P$ such that a parameter of interest $\mu(P)=m$.
A confidence set is a function $C(\cdot)$ of the data $Z$.
The C.S. is said to be $P^*$-a.s. continuous if the boundary of the set $\{m \in C(Z)\}$
has zero probability under $P^*$ for every value of $m$ in the range of $\mu(\cdot)$.
The model $P^*$ is a limit point of $\mmp(m)$ wrt a metric $d(\cdot,\cdot)$ if there exists a  sequence $\{ P_k\}_k \subseteq \mmp(m)$ such that $d(P_k,P^*) \to 0$.
If $\mmp_0$ is dense in $\mmp$, then $co(\mmp_0)$ is dense in $\mmp$ because $\mmp_0 \subseteq co(\mmp_0)$. Similarly, if $P^\star$ is limit point of $\mmp(m)$, then $P^\star$ is also limit point of $co(\mmp(m))$.
Convergence in TV implies convergence in LP. The converse is not generally true. 
See page 10 for sufficient conditions for the converse to hold. 
}
\end{sidewaysfigure}

The first contribution of this paper is to connect the literature on impossible inference and study the relationships between type A and type B impossibility.
Figure \ref{diagram} at the end of this introduction summarizes the literature along with novel relationships derived in this paper.
To the best of our knowledge, impossibility type A dates back to the 1950s.
In a classic paper, \cite{bahadur} show both types of impossibility in the population mean case.
Any test for distinguishing zero mean from non-zero mean distributions has power limited by size; 
and any a.s. bounded confidence interval for the population mean has error probability equal to one.\footnote{The 
impossibility typically arises due to the richness of models in the class of all likely models.
Impossibility does not arise if we restrict the class to only one model, which is the same as pointwise inference.
Uniform inference over a larger class of models is important because the researcher typically does not know all aspects of the model at hand. 
For example, instruments could be weak, and if we incorrectly assume they are always strong, pointwise inference conclusions are quite misleading. }
\cite{bahadur} employ the Total Variation (TV) metric to measure the distance between any two distributions.
We refer to this notion of distance as \emph{strong distance}.
They show that the null set of distributions with a certain mean is dense with respect to (wrt) the TV metric in the set of distributions with all possible means. 

In fact, impossibility type A is very much related to the density of the convex hull of the null set in the set of all likely models wrt the TV metric.
\cite{kraft1955} targets the problem of testing any two sets of distributions and arrives at an important generalization of the theory of \cite{bahadur}.
Kraft's Theorem 5 gives a necessary and sufficient condition for the existence of a test whose  minimum power is strictly greater than its size. 
Such a test exists if, and only if, the minimal TV distance between the convex hulls of the null and alternative sets is bounded away from zero.
Kraft attributes the theorem to Le Cam, and an analogous version of his theorem appears as Theorem 2.1 of \cite{ingster2003}.
\cite{romano2004} demonstrates that the null set being dense in the set of all likely models wrt the TV metric is a sufficient condition for impossibility A.
We derive a corollary of Kraft's Theorem 5 that says that the convex hull of the null set being dense in the set of all likely models wrt the TV metric is a necessary and sufficient condition for impossibility type A.
The null set being dense implies that the convex hull of the null set is dense.
Our corollary connects the literature on impossibility type A wrt the TV metric.

A different branch of the econometrics literature focuses on impossibility type B of confidence sets 
for a given parameter of interest, e.g.  mean or regression slope. 
In the population mean case, \cite{bahadur} arrive at impossibility type B by demonstrating the following fact.
For any mean value $m$, the set of distributions with mean equal to $m$ is dense in the set of all likely models wrt the TV metric.
This is stronger than the sufficient condition for impossibility type B used by \cite{gleser1987}.
\cite{gleser1987} consider classes of models indexed by parameters in a Euclidean space. 
They obtain impossibility type B whenever there exists one distribution $P^*$ such that, for every value of the parameter of interest,
$P^*$ is approximately equal to distributions with that value of the parameter of interest wrt the TV metric.\footnote{\cite{gleser1987} restrict their analysis to
distributions that have parametric density functions wrt the
same sigma-finite measure. Two distributions are indistinguishable if their
density functions are approximately the same pointwise in the data. In their
setting, pointwise approximation in density functions is the same as
approximation in the TV metric. However, pointwise approximation in density
functions is still stronger than convergence in distribution. 
See Proposition 2.29 and Corollary 2.30, \cite{van2000}. 
} 
As with impossibility type A, impossibility type B also holds if the condition of \cite{gleser1987} holds over the convexified space of distributions, which is a weaker sufficient condition.
\cite{donoho1988} also provides type B impossibility for parameters of interest that
are functionals of distributions satisfying a dense graph condition in the TV metric.
One example of such a functional is the derivative of a probability density function (PDF).
Althought it is impossible to obtain a.s. bounded confidence sets, \cite{donoho1988} shows that it is possible to build valid one-sided lower-bounded confidence intervals in some cases.
Similarly, \cite{low1997} demonstrates the impossibility of adaptation gains 
for the length of confidence intervals on linear functionals of non-parametric functions.
Low's lower bound on the expected length of confidence intervals grows to infinity as the class of possible models increases. 

\cite{dufour1997} generalizes the work of \cite{gleser1987} to classes of models indexed by parameters in general metric spaces. 
\cite{dufour1997} also notes that 
impossibility type B implies that tests constructed from a.s. bounded confidence sets fail to control size. 
Unlike all authors mentioned thus far, \cite{dufour1997} relies on a notion of distance much weaker than the TV metric,
which is the notion of distance behind weak convergence or convergence in distribution. 
He obtains impossibility type B whenever there exists one distribution $P^*$ such that, for every value of the parameter of interest,
there exists a sequence of distributions with that value of the parameter of interest that converges in distribution to $P^*$.
The weaker notion of distance restricts the analysis  to confidence sets whose boundary has zero probability under $P^*$.
The L\'{e}vy-Prokhorov (LP) metric is known to metrize weak convergence.
We refer to this notion of distance as \emph{weak distance}.
We demonstrate that the impossibility type B of Dufour also holds after convexifying the space of distributions.

We revisit impossibility type A when distributions are indistinguishable in the LP metric as opposed to the TV metric.
We find that impossibility type A applies to all tests that are a.s. continuous under alternative distributions.
A sufficient condition is that the convex hull of the null set is dense in the set of all likely models wrt the LP metric.
On the one hand, the LP metric does not yield impossibility type A for every test function.
On the other hand, the class of a.s. continuous tests includes the vast majority of tests used in empirical studies.
Convergence in the TV metric always implies convergence in the LP metric. 
The converse is not true, except in more restricted settings. 
For example, if convergence in distribution implies uniform convergence of probability density functions (PDF), then Scheff\'{e}'s Theorem implies convergence in the TV metric (Corollary 2.30 of \cite{van2000}).

The  second contribution of this paper is to note that a weaker notion of distance, such as the LP metric, brings further insights to the problem of impossible inference.
First, it is often easier to prove convergence of models in terms of the weak distance than  in terms of the strong distance.
Application of arguments similar to Portmanteau's theorem immediately yields the LP version of impossible inference in an important class of models in economics that rely on discontinuities. 
Second, the use of the LP metric helps researchers look for tests with non-trivial power. 
If models are indistinguishable wrt the LP metric, but distinguishable wrt the TV metric, we show that a useful test must necessarily be a.s. discontinuous.
Third, the LP metric can be a sensible choice of distance to study hypothesis tests that are robust to small model departures. 
For example, consider the null set of continuous distributions versus the alternative set of discrete distributions with finite support in the rational numbers. 
It is possible to approximate any such discrete distribution by a sequence of continuous distributions in the LP metric.
Hence, it is impossible to powerfully test these sets with a.s. continuous tests.
On the other hand, a positive TV distance between null and alternative leads to a perfect test that rejects the null if observations take rational values.  
Robustness leads us to ask whether observing rational numbers is indeed evidence against the null hypothesis, or simply a matter of rounding or measurement error.
The same problem may arise in reduced-form or structural econometric models, even when the degree of misspecification is small.
Depending on the problem at hand, we may want to look for tests that separate the closure of each hypothesis wrt the LP metric.

The third contribution of this paper is to point out impossible inference in microeconometric models based on discontinuities and macroeconometric models of time series.
Numerous microeconometric analyses identify parameters of interest by relying on natural
discontinuities in the distribution of variables. This is the case of
Regression Discontinuity Designs (RDD), an extremely popular identification
strategy in economics. In RDD, the assignment of individuals into a program
changes discontinuously at a cutoff point in a variable such as age or test
score, as for \cite{hahn2001id} and \cite{imbens2008}. For example, \cite{schmieder2012} study individuals whose duration of unemployment insurance
jumps wrt age. \cite{jacob2004} analyze the effect of students' 
participation in summer school, which changes discontinuously wrt test
scores. Assuming all other characteristics vary smoothly at the cutoff, the
effect of the summer school on future performance is captured by a
discontinuous change in average performance at the cutoff. 
A fundamental assumption
for identification is that performance varies smoothly with test scores, after controlling for summer school.
Models with continuous effects
are well-approximated by models with discontinuous effects. \cite{kamat2015}
uses the TV metric to show that the current practice of tests in RDD suffers
from impossibility type A. We revisit his result using the
LP metric, and we show that impossibility type B
also holds in RDD. 

A Monte Carlo experiment shows that the usual implementation of  Wald tests in RDD, as suggested by \cite{cattaneo2014calonico}, may have size above the desired significance level, even under sensible model restrictions.
We rely on data-generating processes that are consistent with the empirical example of \cite{lee2008}.
Moreover, the simulations show that the Wald test has very little power, even after artificially controlling size.
Slope restrictions on the conditional mean functions do not correct the finite sample failure of the typical Wald test.

In other applications, researchers assume a discontinuous change in
unobserved characteristics of individuals at given points. This is the idea
of bunching, widely exploited in economics. Bunching may occur because of a
discontinuous change in incentives or a natural restriction on variables.
For example, the distribution of reported income may display a non-zero probability 
at points where the income tax rates change, as noted by  
\cite{saez2010}; or, the distribution of average smoking per day 
has a non-zero mass at zero smoking. 
We show that the problem of testing for existence of bunching in a scalar variable 
suffers from type A impossibility for a.s. continuous tests but not for discontinuous tests.

\cite{caetano2015} uses the conditional distribution of variables with bunching and proposes an exogeneity test without instrumental variables.
The key insight is that bunching in the distribution of an outcome variable given a treatment variable constitutes evidence of endogeneity.
For example, consider the problem of determining the effect of smoking on birth weight.
A crucial
assumption is that birth weight varies smoothly with smoking while controlling for all other factors.
Under this assumption, bunching is equivalent to the
observed average birth weight being discontinuous at zero smoking. The
exogeneity test looks for such discontinuity as evidence of endogeneity. Our
point is that models in which birth weight is highly sloped or even
discontinuous based on smoking are indistinguishable from smooth models.
Therefore, we find the exogeneity test has power limited by size. The
current implementation of tests for the size of discontinuity leads to
bounded confidence sets, so it also fails to control size.

In addition to these applications with discontinuities, we verify the existence of impossible inference in a macroeconometrics example where data are continuously distributed. 
We first show that the choice of the weak versus the strong distance connects to the work of Peter J. Huber on robust statistics and leads
us to look at the closure of the set of covariance-stationary time-series processes wrt the LP metric.
This closure contains error-duration models and Compound Poisson models. 
Our theory implies that it is impossible to robustly distinguish these models from covariance-stationary models, even with discontinuous tests.

Our goal with these applications is not to say that valid inference is never possible.
Rather, we point practitioners to the need to restrict either the class of models under consideration or the null hypothesis being tested.
In the RDD case, impossibility vanishes if we restrict the variation of conditional mean functions on either side of the cutoff.
\cite{kamat2015} demonstrates the asymptotic validity of Wald tests under uniform bounds on the derivatives of the conditional mean functions.
\cite{armstrongkolesar2015} derive minimax optimal-length confidence intervals
in the case of a convex class of conditional mean functions, which covers most smoothness or
shape assumptions used in econometrics.
Alternatively, instead of limiting the whole class of models, researchers may consider null hypotheses that restrict other aspects of the model, beyond simply the effect at the threshold.
For example, the null of smooth models with zero effect, or the null of no treatment spillover,
does not suffer from type A impossibility.
We expand this discussion in Section \ref{sec:app:rd} with empirical examples in RDD.

The rest of this paper is divided as follows. Section \ref{sec:theo} sets up
a statistical framework for testing and building confidence sets. 
It presents necessary and sufficient conditions for impossible inference in general non-parametric settings.
Section \ref{sec:rob} connects the LP metric to robust hypothesis testing.
Section \ref{sec:app} gives multiple economic applications where both types of impossibility arise.
Section \ref{sec:simul} presents a Monte Carlo simulation for an empirical application of RDD.
Section \ref{sec:con} concludes.
An appendix contains all formal proofs.
Figure \ref{diagram} (on the next page) summarizes the literature on impossible inference, along with implications of this paper.

\section{Impossible Inference}

\label{sec:theo}

\indent

The researcher has a sample of $n$ observations $Z=(Z_1,\ldots,Z_n)$ that take values in $\m{Z}$, a subset of the Euclidean space $\mmr^{n \times l}$.
The data $Z$ follow a distribution $P$, and the set of all possible distributions considered by the researcher is $\mathbb{P}$.
Every probability distribution $P\in \mathbb{P}$ is defined on the same sample space $\m{Z}$ with Borel
sigma-algebra $\mathcal{B}$. 
It is assumed that all distributions in $\mmp$ are absolutely continuous wrt the same sigma-finite measure.\footnote{Examples include Lebesgue measure for continuous distributions; counting measure for discrete distributions; and sum  of Lebesgue and counting measures for mixed continuous-discrete distributions.} 
We are interested in testing the null
hypothesis $H_{0}:P\in \mathbb{P}_{0}$ versus the alternative hypothesis 
$H_{1}:P\in \mathbb{P}_{1}$ for a partition $\mathbb{P}_{0},\mathbb{P}_{1}$
of $\mathbb{P}$.
 We characterize a hypothesis test by a function of the data 
$\phi : \m{Z} \rightarrow \lbrack
0,1]$. If $\phi $ takes on only the values 0 and 1, the test is said to be 
\emph{non-randomized}, but \emph{randomized} otherwise.
 Given a sample $Z$, we reject
the null $H_{0}$ if the function $\phi (Z)$ equals one, but we fail to
reject $H_{0}$ if $\phi (Z)=0$. If the function $\phi (Z)$ is between 0 and
1, we reject the null with probability $\phi (Z)$ conditional on $Z$. The
unconditional probability of rejecting the null hypothesis under
distribution $P\in \mathbb{P}$ is denoted $\mathbb{E}_{P}[\phi ]$.

The size of the test $\phi $ is $\sup_{P\in \mathbb{P}_{0}}\mathbb{E}_{P}[\phi ]$.
The power of the test under distribution $Q\in \mathbb{P}_{1}$ is given by 
$\mathbb{E}_{Q}[\phi ]$. We say a test $\phi $ has power limited by size when 
$\sup_{Q\in \mathbb{P}_{1}}\mathbb{E}_{Q}[\phi ]\leq \sup_{P\in \mathbb{P}%
_{0}}\mathbb{E}_{P}[\phi ]$. 
Define $co(\mmp')$ to be the convex hull of an arbitrary subset $\mmp' \subseteq \mmp$. That is,
\begin{gather}
co(\mmp') = \left\{ 
P^* : 
~ P^* = \sum_{i=1}^N \alpha_i P_i,
~\text{for some }  N \in \mmn,  P_i \in  \mmp' ~\forall i, 
\right.
\nonumber
\\
\left.
~\alpha_i \in [0,1] ~\forall i, \sum_{i=1}^N \alpha_i =1
\right \}.
\label{eq:convex}
\end{gather}

A small distance between models in $\mathbb{P}_{0}$ and $\mathbb{P}_1$ determines testing impossibility.
There exist various notions of distance to measure the difference between two distributions P and Q. A common choice in
the literature on testing impossibility is the Total Variation (TV) metric $%
d_{TV}(P,Q)$: 
\begin{equation}
d_{TV}(P,Q)=\sup_{B\in \mathcal{B}}\left\vert P(B)-Q(B)\right\vert .
\label{def:tv-metric}
\end{equation}%

Theorem 5 of \cite{kraft1955} says that there exists a test $\phi$ with minimum power strictly greater than
size if, and only if, there exists $\eps>0$ such that $d_{TV}(P,Q) \geq \eps$ for every $P \in co(\mmp_0)$ and $Q\in co(\mmp_1)$.
We restate his theorem below for convenience.
\begin{theorem}\label{theo:kraft}\textbf{(\cite{kraft1955})}
Fix $\eps > 0$. The following statements are equivalent:
\begin{enumerate}
\item[(a)] $ \exists \phi ~:~ 
\inf\limits_{ Q \in \mmp_1 } \mme_Q \phi  \geq \eps + \sup\limits_{ P \in \mmp_0 } \mme_P  \phi$, and

\item[(b)] $\forall P \in co(\mmp_0), ~ \forall Q \in co(\mmp_1), ~ d_{TV}(P,Q) \geq  \eps$.
\end{enumerate}
\end{theorem}

An important implication of Theorem \ref{theo:kraft} for impossible inference is that it gives a necessary and sufficient condition in terms of the convex hull of the null set being dense in the set of all likely models wrt the TV metric.
In other words, 
the convex hull $co(\mathbb{P}_{0})$ is \emph{indistinguishable} from (or dense in) the set of all likely models wrt the TV metric
if, for any $Q\in 
\mathbb{P}_{1}$, there exists a sequence $\{P_{k}\}_{k=1}^{\infty }$ in $co(\mathbb{P}_{0})$ such that $d_{TV}(P_{k},Q)\rightarrow 0$. We demonstrate this fact in the corollary below.
\begin{corollary}\label{coro:kraft_romano}
The following statements are equivalent:
\begin{enumerate}[(a)]
\item for every $Q \in \mmp_1$, there exists a sequence $\{P_k\}_k \subseteq co(\mmp_0)$
such that $d_{TV}(P_k,Q ) \to 0$, and 
\item for every $\phi$ and $Q \in \mmp_1$, 
$ \mme_Q \phi \leq \sup_{P \in \mmp_0 } \mme_P \phi $.
\end{enumerate}
\end{corollary}

The proof of this corollary, as well as all other proofs for the paper, is included in the appendix.
The striking result of \cite{kraft1955} stated in Theorem \ref{theo:kraft}  makes the impossibility type A  found  by \cite{bahadur} and \cite{romano2004} special cases of Corollary \ref{coro:kraft_romano}.
In particular, Theorem 1 of \cite{romano2004} says that Corollary \ref{coro:kraft_romano}-(a) without convexification is a sufficient condition for Corollary \ref{coro:kraft_romano}-(b). 
Notably, \cite{romano2004} finds a positive result for testing population means.
He demonstrates that the t-test uniformly controls  size in large samples with a very weak uniform integrability type of condition, and that the t-test is also asymptotic minimax optimal.

\cite{dufour1997} uses the notion of distance associated with weak convergence to derive impossibility type B.
We say a
sequence $\{P_{k}\}_{k=1}^{\infty }$ converges in distribution to $Q$, if,
for every $B\in \mathcal{B}$ such that $Q(\partial B)=0$, $%
P_{k}(B)\rightarrow Q(B)$. Here, $\partial B$ is the boundary of a Borel set
B, that is, the closure of $B$ minus the interior of $B$. We denote
convergence in distribution by $P_{k}\overset{d}{\rightarrow }Q$. 
Convergence in distribution is equivalent to convergence in the L\'{e}vy-Prokhorov (LP) metric
(\cite{dudley1976}, Theorem 8.3) :
\begin{gather}
d_{LP}(P,Q) = \inf \{\eps > 0 ~:~ P(A) \leq Q(A^\eps) + \eps \text{ for } A \in \m{B} \} 
\label{def:lp-metric}
\\
\text{ where } A^\eps=\{x: \| x-a \|<\eps \text{ for }a \in A \},
\nonumber
\\
\text{ and } \| \cdot \| \text{ is the Euclidean norm on }\mmr^{n \times l}.
\nonumber
\end{gather}

 Convergence of $%
P_{k}$ to $Q$ in the TV metric implies convergence in distribution.
The converse does not hold, in general.\footnote{For example,
a standardized binomial variable converges in distribution to a
standard normal as the number of trials goes to infinity and the probability
of success is fixed. It does not converge in the TV metric because the
distance between these two distributions is always equal to one. In fact,
consider the event equal to the entire real line minus the support of the
binomial distribution. This event has unit probability under the normal
distribution, but zero probability under the binomial distribution.}
It is necessary to restrict the class of distributions in order for convergence in the TV metric to
imply convergence in distribution. 
For example, suppose that $P_k \dto Q$,
that these distributions have common support $[a,b]$ and PDFs $f_{P_k}$, $f_Q$.
Assume further that $f_{P_k}$ converges uniformly over $[a,b]$.
Then, $f_{P_k}$ converges uniformly to $f_Q$ (Theorem 7.17 of \cite{rudin1976principles}).
Convergence of PDFs implies convergence in the TV metric (Scheff\'{e}'s Theorem, see Corollary 2.30 of \cite{van2000}). For a counter-example where these conditions do not hold, consider the bunching example of Section \ref{sec:app:bun}. The null is the set of distributions with a continuously differentiable CDF. The alternative is the set of distributions with a mass point at $x_0$ but
continuously differentiable CDF otherwise. For any CDF $F_Q$ in the alternative,
there exists a sequence of CDFs $F_{P_k}$ in the null that converges pointwise to $F_Q$, so that convergence in distribution holds. Convergence in TV does not hold because $x_0$ has positive probability under $Q$ but zero probability under $P_k$ for every $k$.
It must be the case that the PDFs $f_{P_k}$  do not converge uniformly.
In fact, $F_Q$ has a jump discontinuity at $x_0$, and the derivative of $F_{P_k}$ at $x_0$ grows without limit as $k \to \infty$.


On the one hand, it is true that the zero TV distance provides a necessary and sufficient condition for testing impossibility.
On the other hand, there are examples of models with non-zero TV distance where it seems sensible that no powerful test should exist.
Section \ref{sec:rob} below formalizes this idea, but we start with a simple example for now.
Consider the null set of continuous distributions versus the alternative set of discrete distributions with finite support in the rational numbers. 
It is possible to approximate any such discrete distribution by a sequence of continuous distributions in the LP metric.
We are led to think the data generated by a null model is observationally equivalent to data generated by an alternative model. 
This motivates us to revisit impossibility type A when distributions are indistinguishable in the LP metric.

\begin{assumption}
\label{assu:dense} For every $Q \in \mathbb{P}_1$, there exists
a sequence $\{ P_k \}_{k=1}^\infty$ in $co(\mathbb{P}_0)$ such that $P_k \dto Q$.
In other words, the convex hull $co(\mathbb{P}_{0})$ is indistinguishable from (or dense in) the set of all likely models wrt the LP metric.
\end{assumption}

Assumption \ref{assu:dense} is a sufficient condition for impossibility type A, as described in Theorem \ref{theo:impossible_test}.
\begin{theorem}
\label{theo:impossible_test} 
If Assumption \ref{assu:dense} holds, then any hypothesis test $\phi(Z) $ that is a.s. continuous under any $Q\in\mathbb{P}_{1}$ has power limited by size.
\end{theorem}

\begin{remark}
As noted by \cite{canay2013}, the topology induced by the LP metric is not fine enough to guarantee convergence
of integrals of any test function $\phi$.
Nevertheless, the class of tests that are a.s. continuous under any 
$Q\in \mathbb{P}_{1}$ can be very large.
For example, take a test that
rejects the null when a test statistic is larger than a critical value: $%
\phi (Z)=I\left( \psi \left( Z\right) >c\right) $. This test is
a.s. continuous if the function $\psi $ is continuous and $%
Q\in \mathbb{P}_{1}$ is absolutely continuous wrt the Lebesgue
measure. 
Theorem \ref{theo:impossible_test} only requires a.s. continuity under the alternative $\mmp_1$,
and the null $\mmp_0$ may still contain discrete distributions.
\end{remark}

\begin{remark}
We do not need to restrict Theorem \ref{theo:impossible_test} to the class
of a.s. continuous tests for every case of $\mmp$.
 For example, consider $\mathbb{P}$
to be a subset of the parametric exponential family of distributions with
parameter $\theta $ of finite dimension. Then, for any test $\phi $, the
power function of $\phi $ is continuous in $\theta $, and Theorem \ref{theo:impossible_test} applies under Assumption \ref{assu:dense} (Theorem 2.7.1, \cite{lehmann2005}).
\end{remark}

\begin{remark}
In many instances, Assumption \ref{assu:dense} holds in both directions.
That is, $\mathbb{P}_{1}$ is indistinguishable from $\mathbb{P}_{0}$, and $%
\mathbb{P}_{0}$ is indistinguishable from $\mathbb{P}_{1}$ in the weak
distance. For example, \cite{bahadur} find that any distribution with mean $m $ is well-approximated by distributions with mean $m^{\prime }\neq m$, and
vice-versa. Section \ref{sec:app} finds the same bidirectionality for
models with discontinuities. If Assumption \ref{assu:dense} holds
in both directions, switching the roles of $\mathbb{P}_{0}$ and $\mathbb{P}%
_{1}$ in Theorem \ref{theo:impossible_test} shows that power is equal to
size.
\end{remark}

\bigskip

It is useful to connect our LP version of testing impossibility with the impossibility of
controlling error probability of confidence sets found by \cite{gleser1987}
and \cite{dufour1997}. Define a real-valued function $\mu :\mathbb{P}%
\rightarrow \mathbb{R}$, for example, mean, variance, median, and so on. The set
of distributions $\mathbb{P}$ is implicitly chosen such that $\mu $ is
well-defined. We consider real-valued functions for simplicity, and results
for $\mu $ with more general ranges are straightforward to obtain. The range
of $\mu $ is $\mu (\mathbb{P})$. Suppose we are interested in a confidence
set for $\mu (P)$ when the true model is $P\in \mathbb{P}$. A confidence set
takes the form of a function $C(Z)$. For a model $P\in \mathbb{P}$,
 the \emph{coverage probability} of $C(Z)$ is given by $P \left[ \mu (P) \in C(Z) \right] $.
 The confidence region $C(Z)$ has \emph{confidence level} $1-\alpha $ (i.e. error probability $\alpha $) 
 if $C(Z)$ contains $\mu (P)$ with probability at least $1-\alpha $: 
\begin{equation}
\inf_{P\in \mmp} {P}\left[ \mu (P)\in C(Z)\right] =1-\alpha.
\end{equation}

For any value $m \in \mu(\mmp)$, we define the subset $\mmp(m)$ by
\begin{gather}
\mmp(m)=\{P\in \mmp~:~\mu (P)=m\}.
\label{eq:pm}
\end{gather}

Impossibility type B says that confidence sets that are a.s. bounded under some distributions in $\mmp$
have zero confidence level. The next assumption gives a sufficient condition for impossibility type B in terms of the LP metric.
\begin{assumption}
\label{assu:limpt}
There exists a distribution $P^*$ (not necessarily in $\mmp$) 
such that for every $m \in \mu(\mmp)$ there exists a sequence $\{P_k \}_k$
in $co(\mmp(m))$ such that $P_k \dto P^*$.
\end{assumption}
 
If Assumption \ref{assu:dense} 
holds with $\mmp_{0}=\mmp(m)$ for every $m\in\mu(\mmp)$, then 
Assumption \ref{assu:limpt} holds.
 In fact, if $\mmp(m)$ is dense in $\mmp$ for every $m$,
then Assumption \ref{assu:limpt} is satisfied for $P^*=Q$ for any $Q \in \mmp_1$.
Some models satisfy Assumption \ref{assu:dense} with $\mmp_0=\mmp(m)$ for every $m\in\mu(\mmp)$ and suffer from both types of impossibility.
Examples of this case include the problem of testing the mean (\cite{bahadur}),
or the problem of testing the size of the discontinuity in RDD (Section \ref{sec:app:rd}).
Nevertheless, some other models satisfy Assumption \ref{assu:limpt} but not Assumption \ref{assu:dense} for every $m$.
These models suffer from impossibility type B.
Examples include the problem of ratio of regression parameters (\cite{gleser1987}), and the problem of weak instruments (\cite{dufour1997}).

The next theorem encapsulates the impossibility of controlling coverage probabilities found by 
\cite{gleser1987} and \cite{dufour1997}. 
It differs from \cite{gleser1987} because Assumption \ref{assu:limpt} uses the LP distance.
It differs slightly from \cite{dufour1997} because
Assumption \ref{assu:limpt} is stated in terms of the convex hull of $\mmp(m)$ rather than simply 
$\mmp(m)$.
 \begin{theorem}
\label{theo:ci} 
 Suppose Assumption \ref{assu:limpt} holds with $P^*$.
 Assume the confidence set $C(Z)$ 
 has confidence
level $1-\alpha $, and $P^*(\partial\{m \in C(Z) \})=0$ for every $m\in \mu(\mmp)$.
Then,
\begin{gather}
\forall m \in \mu(\mmp) ~:~ P^* \left[ m \in C(Z) \right] \geq 1-\alpha.
\label{theo:ci:eq1}
\end{gather}
For a set $A \subset \mmr$, 
define $U[A] = \sup\{c: c\in A \}$, $L[A] =  \inf\{c: c\in A \}$,
and $D[A] = U[A] - L[A]$.
Assume $\{ U[C(Z)] \geq x \}$, $\{ L[C(Z)] \leq -x \}$,
and $\{ D[C(Z)] \geq x \}$ are measurable events for every $x \in [0,\infty]$.
If $D[\mu(\mmp)] = \infty$, then 
\begin{gather}
P^* \left[ D[ C(Z) ] =\infty  \right] \geq 1-  \alpha.
\label{theo:ci:eq2}
\end{gather}
In addition, if  $P^*\left[ \partial \{ D[C(Z)] =\infty   \} \right]=0$, then  
\begin{gather}
\forall \eps >0  ~:~ \sup\limits_{P \in B_{\eps}(P^*) \cap \mmp} 
P\left[ D[C(Z)] = \infty \right] \geq 1- \alpha
\label{theo:ci:eq3}
\end{gather}
where $B_{\eps}(P^*)= \{P  ~:~ d_{LP}(P,P^*)<\eps \}$. 
 \end{theorem}

\begin{remark}
Part (\ref{theo:ci:eq3}) above implies the following.
If $1-\alpha>0$, then the
confidence set $C(Z)$ is unbounded with strictly positive probability for
some  $P\in \mathbb{P}$.
Alternatively, the contrapositive of part (\ref{theo:ci:eq3}) 
says the following. 
Any confidence
set that is a.s. bounded under distributions in $\mmp$ in a neighborhood of $P^*$
has $1-\alpha=0$ confidence level. 
\end{remark}
\begin{remark}
It is possible to obtain a slightly more general version of Theorem \ref{theo:ci} using Assumption \ref{assu:limpt} stated in terms of the TV metric as opposed to the LP metric.
In that case, Theorem \ref{theo:ci} would be true for confidence sets that do not necessarily satisfy $P^*(\partial\{m \in C(Z) \})=0$ and
$P^*\left[ \partial \{ D[C(Z)] =\infty   \} \right]=0$.
\end{remark}

A common way to obtain confidence sets is to invert hypothesis tests.
The function $C(Z)$ is constructed by inverting a test in the following
manner. For a given $m\in \mu ({\mathbb{P}})$, define 
$\mmp_{0,m} = \mmp(m)$ and $\mmp_{1,m}=\mmp \setminus \mmp(m)$,
where 
$A \setminus B$ denotes the remainder of set $A$ after we remove the intersection of set $B$ with set $A$.
If $\phi_m(Z)$ is a non-randomized test\footnote{In case the test $\phi_m(Z)$ is randomized, 
use $\ti{\phi}_m(Z) = \mmi\{U \leq \phi_m(Z) \}$, where $U$ is uniformly distributed in $[0,1]$ and  independent of $Z$.} for $\mmp_{0,m}$ vs $\mmp_{1,m}$, 
then
\begin{equation}
C(Z)=\{m\in \mu(\mmp)~:~\phi _{m}(Z)=0\}.  
\label{eq:theo:ci_test}
\end{equation}

For every $m\in \mu( \mmp )$, the test $\phi_m(Z)$ has size $\alpha(m) = \sup_{P \in \mmp_{0,m}}\mme_{P}\left[\phi_m(Z)\right]$.
The confidence level of $C(Z)$ is equal to one minus the supremum of $\alpha(m)$ over $m \in \mu(\mmp)$.
The proof of this claim is found in Lemma \ref{lem:covsize} in the appendix.

Theorem \ref{theo:ci} along with Lemma \ref{lem:covsize} imply that tests that invert into a.s. bounded confidence sets fail to control size.
\begin{corollary}
\label{coro:ci}
Suppose Assumption \ref{assu:limpt} holds,
and $\mu (\mathbb{P})$ is unbounded. 
Let the confidence set $C(Z)$ be constructed from tests $\phi_m(Z)$, as in 
Equation (\ref{eq:theo:ci_test}).
Assume $C(Z)$ has confidence level $1-\alpha $ and satisfies the assumptions of Theorem \ref{theo:ci}.
If $C(Z)$ is a.s. bounded under distributions in $\mmp$ in a neighborhood of $P^*$,
then $\alpha=1$.
Consequently, for every $\eps>0$, there exists $m_{\eps} \in \mu(\mmp)$
such that $\sup\limits_{P \in \mmp_{0,m_{\eps}}} \mme_P \phi_{m_{\eps}} > 1 -\eps$.
\end{corollary}

\begin{remark}
\cite{moreira2003} provides numerical evidence that Wald tests can have
large null rejection probabilities for the null of no causal effect ($m=0$)
in the simultaneous equations model. To show that Wald tests have null
rejection probabilities arbitrarily close to one, the hypothesized value $m$
for the null would need to change as well. He also suggests replacing the
critical value by a critical value function of the data. This critical value
function depends on the hypothesized value $m$. Our theory shows that this
critical value function is unbounded if we change $m$ freely.
\end{remark}

\section{Weak Convergence and Robustness}
\label{sec:rob}

\indent

This section presents further motivation for using the LP metric to study impossible inference.
It relates the weak topology induced by the LP metric
to the theory developed by Peter J. Huber, 
the most prominent researcher
in the area of robust statistics. 
We refer the reader to \cite{huber2009} for more details.
We start this section with a discussion of robust statistical procedures.
An example of impossible robust hypothesis testing in time-series models appears in Section \ref{sec:app:ts}.

Several statistical procedures are susceptible to small model departures.
This perception has led researchers to propose alternative procedures that are less sensitive to the break down of usual assumptions.
Huber studies different ways of defining a set of model departures $\mmp_{\eps}$.
One possibility is to assume that the actual distribution of the data is a mixture of a distribution in $\mmp$ with a distribution from a more general set of models $\mmm$.
In other words, $\mmp$ may be contaminated with probability $\eps $:
\begin{equation}
\mathbb{P}_{\eps }=\left\{ H\in \mathbb{M};~ \exists F\in \mathbb{P}\text{
and }\exists G\in \mathbb{M};~H=\left( 1-\eps \right) F+\eps
~G\right\} ,  \label{(eq epsilon contaminated P)}
\end{equation}%
where $\mathbb{M}$ is  larger than the original $\mathbb{P}$. Estimators or tests are said to be robust if they have minimax properties over the set of model departures 
$\mmp_{\eps}$.
To highlight the importance of robust procedures, we briefly discuss two examples. 

The first example of a robust procedure involves point estimation.
The researcher has a sample of $n$ iid observations $Z_{i}\in \mathbb{R}^{l}$, $i=1,\ldots ,n$.
The set of joint probability distributions $\mathbb{P}$ is indexed by a parameter $\theta $ and admits marginal densities $p\left( Z_{i};\theta \right) $
wrt the same dominating measure (e.g. Lebesgue).
The maximum likelihood estimator (MLE) then minimizes
\begin{equation*}
\sum\nolimits_{i=1}^{n}-\ln p\left( Z_{i};\theta \right) .
\end{equation*}
This estimator $\ha{\theta}$ solves
\begin{equation*}
\sum\nolimits_{i=1}^{n}-\frac{\partial p\left( Z_{i}; \ha{\theta} \right) }{\partial \theta }.\frac{1}{p\left( Z_{i};\ha{\theta} \right) }=0.
\end{equation*}
Under the usual regularity conditions, $\ha{\theta}$ is consistent, asymptotically normal, and efficient within the class of regular estimators.

A common choice for $\mmm$ is the set of distributions with symmetric, thick-tailed densities.
It is well-known that optimal procedures derived under Gaussian distributions
(sample drawn from $\mathbb{P}$) break down if there is a probability $\eps $
 of observing outliers (sample drawn from $\mathbb{M}$).
\cite{huber1964} suggests M-estimators. To give a specific example of a robust M-estimator, consider the regression model
\begin{equation*}
Y_{i}=X_{i}^{\prime }\theta +U_{i}\text{,}
\end{equation*}%
where we observe $Z_{i}=\left( Y_{i},X_{i}\right) $ but do not observe the
zero-mean normal errors $U_{i}$. The MLE $\ha{\theta}$ minimizes 
\begin{equation*}
\sum\nolimits_{i=1}^{n} \left( Y_{i}-X_{i}^{\prime }\theta\right)^2 \text{,}
\end{equation*}%
and satisfies 
\begin{equation*}
\sum\nolimits_{i=1}^{n}X_{i} ~ \left( Y_{i}-X_{i}^{\prime } \ha{\theta} \right) =0.
\end{equation*}

More generally, a M-estimator $\ha{\theta}$ minimizes 
\begin{equation*}
\sum\nolimits_{i=1}^{n}\rho\left( Y_{i}-X_{i}^{\prime }\theta \right) 
\text{,}
\end{equation*}
 and satisfies 
\begin{equation*}
\sum\nolimits_{i=1}^{n} X_{i} ~\psi\left( Y_{i}-X_{i}^{\prime }\ha{\theta}\right) =0
\end{equation*}%
for choices of functions $\rho$ and $\psi$.
In the MLE case above, $\rho(u)=u^2$ and $\psi(u)=u$.

An M-estimator $\ha{\theta}$ is said to be asymptotically minimax optimal among a class of estimators if it minimizes the maximal asymptotic variance over distributions in 
$\mmp_{\eps}$. 
The M-estimator associated with the functions
\begin{equation}
\rho _{k}\left( u\right) =\left\{ 
\begin{array}{cc}
u^{2}/2 & \text{if }\left\vert u\right\vert \leq k \\ 
k\left\vert u\right\vert -u^{2}/2 & \text{if }\left\vert u\right\vert >k%
\end{array}%
\right. \text{ and }\psi _{k}\left( u\right) =\max \left\{ -k,\min \left(
k,u\right) \right\}   \label{(eq pho_k)}
\end{equation}%
is known to be asymptotically minimax optimal for model contamination.
The constant $k$ depends on the deviations $\eps $ in (\ref{(eq epsilon
contaminated P)}). 
 As $\eps \rightarrow 0$, the truncation parameter $k\rightarrow \infty $.
 As the model departure is small, the M-estimator approaches the MLE estimator. 
 If $\eps \rightarrow 1$, the parameter $k\rightarrow 0$. 
 As the contamination is arbitrarily large, the M-estimator approaches the least absolute deviation (LAD) estimator.

The second example of a robust procedure is in hypothesis testing.
Consider the problem  of  testing a simple null $P_{0}$ against a simple alternative $P_{1}$.
 Assume both $P_{0}$ and $P_{1}$ have densities $p_{0}$ and $p_{1}$ wrt the Lebesgue measure.
 For a sample $X=\left(X_{1},...,X_{n}\right) $, the likelihood ratio (LR) test rejects the null if
and only if
\begin{equation*}
\prod\nolimits_{i=1}^{n}\frac{p_{1}\left( X_i \right) }{p_{0}\left( X_i \right) }
>c_{\alpha },
\end{equation*}%
where $c_{\alpha }$ is the $1-\alpha $ quantile of the distribution of the
left-hand side under the null. 
The Neyman-Pearson Lemma asserts that the LR test is optimal, as it maximizes power within the class of tests with correct size $\alpha $.

Similar to model departures in the point-estimation example above, we consider the possibility that the null and alternative hypotheses are misspecified.
The $\eps $-contaminated null and alternatives are
\begin{equation}
\mathbb{P}_{i,\eps }=\left\{ H\in \mathbb{M};\exists F\in \mathbb{P}_{i}\text{ and }\exists G\in \mathbb{M};H=\left( 1-\eps \right) F+\eps ~G\right\} ,  \label{(eq null and alternative contaminated sets)}
\end{equation}%
for $i=0,1$. The new sets $\mmp_{0,\eps }$ and $\mmp_{1,\eps }$ allow for local departures for arbitrary distributions in $\mmm$.
A minimax optimal hypothesis test maximizes the minimal power over $\mmp_{1,\eps }$ within the class of tests with correct size over $\mmp_{0,\eps }$.

\cite{huber1965} shows that the minimax test to these model departures rejects the null if and only if 
\begin{equation*}
\prod\nolimits_{i=1}^{n}\pi _{k}\left( \frac{p_{1}\left( x\right) }{%
p_{0}\left( x\right) }\right) >c_{\alpha },
\end{equation*}%
where
\begin{equation*}
\pi _{k}\left( w\right) =\max \left\{ k_{1},\min \left( k_{2},w\right)
\right\} ,
\end{equation*}%
for constants $k=\left( k_{1},k_{2}\right) $ that depend on the size of the departure $\eps $.
As $\eps \rightarrow 0$, the constant $k_{1}$ approaches zero, and $k_{2}$ diverges to infinity. 
Hence, as the departure decreases,  the robust test approaches the usual LR test.

In the two examples of robust procedures given above, arbitrarily small model departures ($\eps \rightarrow 0$) do not affect the solution to the minimax problem.
That is, as $\eps$ approaches zero, the robust estimator converges to the MLE, and the robust test converges to the LR test.
These limiting solutions remain the same, even if we ignore model departures ($\eps=0$).
These inference procedures target a parameter which is a functional of the underlying distribution $P \in \mmp_{\eps}$.
These functionals vary smoothly wrt $\eps$ as $\eps \to 0$.
Robustness is associated with smoothness of the functional, but such smoothness may not always occur in other settings.

Our work on impossible inference and the different metrics connects to Huber's work on robustness when we look at the following definition of model departure.
For a metric space $(\mmm,d)$, the set of model departures is defined as an $\eps$-neighborhood of $\mmp$:
\begin{equation*}
\mathbb{P}_{\eps }=\left\{ H\in \mathbb{M}\mathcal{;\exists }F\in \mathbb{P}\text{ s.t. }d\left( F,H\right) \leq \eps \right\} .
\end{equation*}

The set $\mathbb{P}_{\eps}$ is closed.\footnote{In 
fact, take an arbitrary convergent sequence $H_{n}\rightarrow H$, such that $H_n \in \mmp_{\eps} ~\forall n$. 
To show $H \in \mmp_{\eps}$, pick an arbitrary $F \in \mmp$.
It is true that $d(H_n,F)\leq \eps ~\forall n$.
Therefore, $d\left( F, H \right) \leq d\left( F,H_{n}\right) +d\left( H_{n},H\right) \leq \eps +d\left( H_{n},H\right)$.
Taking the limit as $n \to \infty$ gives $d\left( F, H \right) \leq \eps$.
}
The set $\bigcap _{\eps >0}\mathbb{P}_{\eps }$ is also closed and coincides
with $\overline{\mathbb{P}}$, the closure of $\mathbb{P}$. 
Hence, the set $\overline{\mathbb{P}}$ is the minimal set of the Huber-type model departures $\mathbb{P}_{\eps }$ containing $\mathbb{P}$.

The minimal set of model departures crucially depends on a choice for the metric $d$.
Aside from the L\'{e}vy-Prokhorov (LP) and the Total Variation (TV) metrics, there are many choices of metrics for spaces of probability measures:  Kolmogorov, Hellinger, and Wasserstein, among others. \cite{gibbs2002} provide a review. 
Which metric shall we choose?
 The choice of the metric on the space of models $\mmm$ induces a topology $\m{V}$ on that space.
Parameters of interest are functionals $\mu:\left( \mathbb{M},\mathcal{V}\right) \rightarrow \left( \mathbb{R},\mathcal{U}\right)$
where $\m{U}$ is the topology on the $\mmr$ space.
Robustness is about the continuity of the functional $\mu$, which crucially depends on the choices of topologies $\mathcal{V}$ and $\m{U}$.
For the real line, it seems reasonable to work with the smallest topology involving all open sets of the form $\left( a,b\right)$. 
However, there are many choices of topologies for the set of measures $\mmm$.

As the topology becomes finer on the domain of $\mu$, the set of continuous functionals $\mu$ grows larger. 
Let us consider a simple example to illustrate this point.
Take two topological spaces, $\left(\mathbb{R},\mathcal{V}\right) $ and $\left( \mathbb{R},\mathcal{U}\right) $,
and a function $\mathbb{\psi }\left( x\right) =x$. Continuity of this simple function requires $\psi^{-1}(U) \in \m{V}$
for every open set $U \in \m{U}$.
Take $U=\left( 0,1\right) $, then $\mathbb{\psi }^{-1}\left( U\right) =\left( 0,1\right) $. 
If we choose the coarsest topology $\mathcal{V}\mathbb{=}\left\{ \emptyset ,\mathbb{R}\right\} $, then even this simple function is not continuous.
It seems reasonable to require all linear functions to be continuous. 
If the topology $\mathcal{V}$ is generated by all open sets of the form $\left( a,b\right) $, then all linear functions are continuous.  
Of course, other non-linear functions may be continuous as well; e.g., $\mathbb{\psi }\left( x\right) =x^{2}$.
This example makes the point that continuity depends heavily on the topology $\mathcal{V}$ associated to the domain of the function. 
If a function is continuous for a topology $\mathcal{V}$, then it is also continuous for a finer topology. 
This relates back to the discussion of robustness as continuity of a functional $\mathbb{\mu }$.

If we choose a fine topology, then many statistical procedures will be deemed robust, because many functionals will be continuous.
If we choose a coarse topology, then fewer statistical procedures will be robust. 
However, if a statistical procedure is robust in the coarse topology, it is also robust in the fine topology.
Among the commonly-used notions of distance in measure spaces, the notion of distance behind weak convergence or convergence in distribution induces the coarsest topology.
The LP metric is a notion of distance that metrizes weak convergence (\cite{dudley1976}, Theorem 8.3). 
To be conservative, if we were to choose one metric, we would choose one that metrizes weak convergence.
After all, if a functional is continuous wrt the topology induced by the LP metric, it is also continuous wrt the stronger topologies induced by the Kolmogorov or TV metrics.

Another question is whether we should be stricter with robustness and look for an even weaker topology than the weak topology induced by the LP metric. 
In perfect analogy to the real line example, the weak topology is the coarsest topology that guarantees continuity for all
functionals of the form
\begin{equation}
\mu \left( P \right) =\int g~dP,
\label{eq:robust:1}
\end{equation}
for $g$ bounded and continuous.
It seems reasonable, after all, to require $\mu$ to be continuous when $g$ is a bounded and continuous function.
If we choose a weaker topology, then not even $\mu$ of this form will be continuous. 

In hypothesis testing, robustness over a minimal set of model departures motivates testing  $\overline{\mathbb{P}}_{0}$ against $\overline{\mathbb{P}}_{1}$
instead of testing $\mathbb{P}_{0}$ against $\mathbb{P}_{1}$.
Allowing for robustified hypotheses $\overline{\mathbb{P}}_{0}$ and $\overline{\mathbb{P}}_{1}$ potentially protects us against numerical approximation errors,
misspecified models, measurement errors, and optimization frictions, among other deviations from the set of models we are testing.
Robustness of inference procedures for $\mu$ that are as simple as (\ref{eq:robust:1}) requires a topology no weaker than the topology induced by the LP metric.
Therefore, we use the LP metric to define the closure of a set for robust hypothesis testing. 
We give two simple examples to strengthen the argument of why the LP metric may be a sensible choice.

The first example of robust hypothesis testing using the LP metric compares extremely simple discrete distributions under both null and alternative hypotheses.
Take $X$ to be a Bernoulli random variable and $X_{n}=X+1/\left( 1+n\right) $ for $n \in \mmn$.
Let $P^X$ denote the distribution of $X$.
The minimal TV distance between $\mathbb{P}_{0}=\left\{P^{X_{n}}\text{ for } n \in \mmn \right\} $
and 
$\mathbb{P}_{1}\mathcal{=}\left\{ P^{X}\right\} $
is equal to one.
 According to Theorem 5 of \cite{kraft1955}, there exists a test for $\mmp_0$ vs $\mmp_1$ with non-trivial power.
For example, define a test which rejects the null if we observe the values $0$ or $1$, but fails to reject the null otherwise.
This test has size equal to zero and power equal to one.
Should we take the values 0 and 1 as evidence against the null? 
Or should we think instead that the null could have led to those same values,  for all practical purposes?
In this example, we note that $\mathbb{P}_{1}\mathcal{\subset }\overline{\mathbb{P}}_{0}$, where the closure is defined wrt the LP metric.\footnote{In fact, 
$F_{n}\left( x\right) =P\left( X_{n}\leq x\right) =P\left( X\leq x-1/n\right)$, 
$F_{n}\left( x\right) \rightarrow P\left( X<x\right) =F\left( x^{-}\right)$,
$F(x^{-}) \neq F(x) ~\Leftrightarrow~ x\in \{0,1\}$
where $ \{0,1\}$ are the only discontinuity points of $F$, so $P^X$ is a limit point of $\mmp_0$.} 
Hence, the minimal TV distance between $\overline{\mathbb{P}}_{0}$ and $\mathbb{P}_{1}$ is equal to zero.
After we robustify the null set  to $\overline{\mathbb{P}}_{0}$, it becomes impossible to find any test with power greater than size 
(Corollary \ref{coro:kraft_romano}).
However, if we define the closure of $\mmp_0$ wrt the TV metric, say $\overline{\mathbb{P}}_{0}^{TV}$, the minimal TV distance between  $\overline{\mathbb{P}}_{0}^{TV}$ and 
$\mathbb{P}_{1}$ is non-zero, which means it is still possible to powerfully distinguish these sets.

The second example of robust hypothesis testing using the LP metric uses the multinomial approximation to continuous distributions.
Take $\mathbb{P}_{0}$ as the collection of multinomial distributions, with each support being a finite subset of rational numbers.
Let $\mathbb{P}_{1}$ be the set of continuous distributions. 
The minimal TV distance between $\mathbb{P}_{0}$ and $\mathbb{P}_{1}$ is one, and it is possible to powerfully distinguish these sets.
Robustness leads us to ask whether observing rational numbers is indeed evidence of the null hypothesis, or simply a matter of rounding  or measurement error.
The closure $\overline{\mathbb{P}}_{0}$ wrt the LP metric contains continuous distributions, 
and the minimal TV distance between $\mathbb{P}_{1}$ and $\overline{\mathbb{P}}_{0}$ is zero.
After we robustify the null set to $\overline{\mmp}_0$, it becomes impossible to powerfully test these hypotheses.

Both in the Bernoulli and multinomial examples, it becomes clear that the LP closure of the null set robustifies the testing procedure.
The next step is to check the TV distance between the robustified null and alternative sets as a way to search for robust tests with non-trivial power.
The use of the TV metric in the second step is justified by a corollary of Theorem 5 of \cite{kraft1955}.
Corollary \ref{coro:kraft_romano}  
demonstrates that a necessary and sufficient condition for the existence of tests with non-trivial power is that the null set is not dense in the set of all distributions wrt the TV metric.

\section{Applications}
\label{sec:app}

\indent

In this section, we apply our theory to multiple economic examples.
The first three examples are of models with discontinuities: RDD, bunching in a scalar variable, and exogeneity tests based on bunching.
In these settings, the proof of the LP version of impossible inference follows arguments similar to Portmanteau's Theorem.
That is, the indicator functions are approximately the same as the steep continuous functions using the weak distance.
The problem of testing for the existence of bunching in a scalar variable differs from the other applications with discontinuities because
there exists a discontinuous powerful test.
A fourth example is in time series; it connects the LP version of impossible inference to Huber's work on robust statistics.
This connection leads to the conclusion that it is impossible to powerfully discriminate error-duration or Compound Poisson models from covariance-stationary models.

\subsection{Regression Discontinuity and Kink Designs}

\label{sec:app:rd}

\indent

The first example is the Regression Discontinuity Design (RDD), first
formalized by \cite{hahn2001id} (HTV01). RDD has had an enormous impact in
applied research in various fields of economics. Applications of RDD started
gaining popularity in economics in the 1990s. Influential papers include those of \cite{black1999better}, who studies the effect of quality of school districts on
house prices, where quality changes discontinuously across district
boundaries; \cite{angrist1999}, who measure the effect of class sizes on
academic performance, where size varies discontinuously with enrollment; and 
\cite{lee2008}, who analyzes US House of Representatives elections and incumbency, where
election victory is discontinuous on the share of votes.

Recent theoretical contributions include the study of rate optimality of RDD
estimators by \cite{porter2003rd} and the data-driven optimal bandwidth
rules by \cite{imbens2012optimal} and \cite{cattaneo2014calonico}. RDD
identifies causal effects local to a cutoff value; 
several authors
develop conditions for extrapolating local effects farther away from the cutoff.
These include estimation of derivatives of the treatment effect at the
cutoff by \cite{dong2014} and \cite{donglewbel2015}; tests for homogeneity
of treatment effects in fuzzy RDD by \cite{bertanha_imbens_2019}; and
estimation of average treatment effects in RDD with variation in cutoff
values by \cite{bertanha2019}. All these theoretical contributions rely on
point identification and inference, and they are subject to both types
of impossibility. The current practice of testing and building
confidence intervals relies on Wald test statistics $(t(Z)-m)/s(Z)$, where $%
t(Z)$ and $s(Z)$ are a.s. continuous and bounded in the data. For a choice
of critical value z, hypothesis tests $\phi(Z) = \mathbb{I}\{ |(t(Z)-m)/s(Z)
| > z \}$ are a.s. continuous when the data is continuously distributed.
Confidence intervals $C(Z)=\{t(Z)-s(Z)z \leq m \leq t(Z)+s(Z)z\}$ have a.s.
bounded length $2 s(Z) z$.

The setup of RDD follows the potential outcome framework.
 For each individual $i=1,\ldots ,n$, define
four primitive random variables $D_{i},X_{i},Y_{i}(0),Y_{i}(1)$. These
variables are independent and identically distributed. The variable $D_{i}$
takes values in $\{0,1\}$ and indicates treatment status. The real-valued
variables $Y_{i}(0)$ and $Y_{i}(1)$ denote the potential outcomes,
respectively, if untreated and treated. Finally, the forcing variable $X_{i}$
represents a real-valued characteristic of the individual that is not
affected by the treatment. The forcing variable has a continuous PDF $f(x)$ with interval support equal to $\mathbb{X}$. The
econometrician observes $X_{i}$, $D_{i}$, and only one of the two potential
outcomes for each individual: $Y_{i}=D_{i}Y_{i}(1)+(1-D_{i})Y_{i}(0)$. For
simplicity, we consider the sharp RDD case, but it is straightforward to
generalize our results to the fuzzy case. In the sharp case, agents receive
the treatment if, and only if, the forcing variable is greater than or equal
to a fixed policy cutoff $c$ in the interior of support $\mathbb{X}$.
Hence, $D_{i}=\mathbb{I}\{X_{i}\geq c\}$, where $\mathbb{I}\{\cdot \}$
denotes the indicator function.

We focus on average treatment effects. In RDD settings, identification of
average effects is typically obtained only at the cutoff value after
assuming continuity of average potential outcomes, conditional on the forcing
variable. In other words, we assume that $\mathbb{E}[Y_{i}(0)|X_{i}=x]$ and $%
\mathbb{E}[Y_{i}(1)|X_{i}=x]$ are bounded continuous functions of $x$. HTV01
show that this leads to identification of the parameter of interest: 
\begin{equation}
m=\mathbb{E}\left[ Y_{i}(1)-Y_{i}(0)|X_{i}=c\right] =\lim\limits_{x%
\downarrow c}\mathbb{E}\left[ Y_{i}|X_{i}=x\right] -\lim\limits_{x\uparrow c}%
\mathbb{E}\left[ Y_{i}|X_{i}=x\right] .  \label{eq:app:rdd:def_m}
\end{equation}

Let $\mathcal{G}$ denote the space of all functions $g: \mathbb{X} \to \mathbb{R}$ that are bounded, and that are infinitely many times continuously differentiable in every $x \in 
\mathbb{X} \setminus \{c\}$.
The notation $\mathbb{X} \setminus \{c\}$ represents the set with every point of $\mathbb{X}$ except for $c$.
Continuity of functions in $\m{G}$ suffices to show impossible inference in this section.
Nevertheless, non-parametric estimators of the size of the discontinuity $m$ typically assume that functions in $g$ are continuously differentiable of first or second order.
We impose that functions in $\m{G}$ are continuously differentiable of infinite order, to demonstrate that both types of impossibility hold even in this more restricted class of functions.
The size of the discontinuity $m$ at the cutoff may take any value in $\mmr$. 

Each
individual pair of variables $Z_i=(X_i , Y_i)$ is iid as $P$. The family of all
possible models for $P$ is denoted as 
\begin{equation}  \label{eq:app:rdd:familyp}
\mathbb{P} = \{P: (X_i, Y_i) \sim P,~ \exists g \in \mathcal{G} \text{ s.t. } 
\mathbb{E}_P[Y_i|X_i=x] = g(x) \}.
\end{equation}

The local average causal effect is the function of the distribution of the
data $P\in \mathbb{P}$ given by (\ref{eq:app:rdd:def_m}), provided the
identification assumptions of HTV01 hold. The parameter $m$ of the size of the
discontinuity is weakly identified in the set of possible true models $%
\mathbb{P}$. 
Intuitively, any conditional mean function $\mathbb{E}[Y_{i}|X_{i}=x]$ that is continuous except for a jump discontinuity at $x=c$ is well-approximated by a
sequence of continuous conditional mean functions. The reasoning behind this approximation is similar to the proof of part of Portmanteau's
theorem (Theorem 25.8, \cite{billingsley2008}). It is known that, if $\mathbb{E}[f(X_{n})]\rightarrow \mathbb{E}[f(X)]$ for every bounded function $f$ that is a.s. continuous under the
distribution of $X$, then $X_{n}\overset{d}{\rightarrow }X$. The proof of Corollary \ref{coro:rdd} uses an infinitely continuously differentiable function $f$ that is approximately equal to an indicator function. 


\begin{corollary}
\label{coro:rdd} Assumption \ref{assu:dense} is satisfied for $\mmp_{0,m}$ $\forall m\in 
\mathbb{R}$, and Theorems \ref{theo:impossible_test} and \ref{theo:ci} apply
to the RDD case. Namely, (i) a.s. continuous tests $\phi _{m}(Z)$ for the
value of the discontinuity $m$ have power limited by size; and (ii) confidence
sets for the value of the discontinuity $m$ and with finite expected length
have zero confidence level.
\end{corollary}

\begin{remark}
Corollary \ref{coro:rdd} also applies to quantile treatment effects by
simply changing the definition of the functional $\mu(\mathbb{P})$ to be the
difference in side limits of a conditional $\tau$-th quantile $%
Q_{\tau}(Y_i|X_i=x)$ at $x=c$.
This contrasts with the problem of testing unconditional quantiles, which does not suffer from impossible inference.
See \cite{lehman1975}, \cite{tibshirani1988},
and \cite{coudin2009}.
\end{remark}

\begin{remark}
In the fuzzy RDD case, the treatment effect is equal to the discontinuity in 
$\mathbb{E}[Y_{i}|X_{i}]$ at $X_{i}=c$ divided by the discontinuity in $%
\mathbb{E}[D_{i}|X_{i}]$ at $X_{i}=c$. Corollary \ref{coro:rdd} applies to
both of these conditional mean functions, and it leads to impossible
inference in the fuzzy RDD case as well. \cite{feir2016} study weak
identification in fuzzy RDD and propose a robust testing procedure.
In contrast to \cite{kamat2015} and to this paper, their source of weak
identification comes from an arbitrarily small discontinuity in $\mathbb{E}%
[D_{i}|X_{i}]$ at $X_{i}=c$.
\end{remark}

The most common inference procedures currently in use in applied research
with RDD rely on Wald tests that are a.s. continuous in the data and produce
confidence intervals of finite expected length. See \cite{imbens2012optimal}
and \cite{cattaneo2014calonico} for the most commonly-used inference procedures.
Corollary \ref{coro:rdd} implies that it is impossible to control size of these
tests and coverage of these confidence intervals.

Ours is not the first paper to show impossible inference in the RDD case. 
\cite{kamat2015} demonstrates that models with a
discontinuity are similar to models without a discontinuity in the TV
metric. He applies the testing impossibility of \cite{romano2004} and finds
that tests have power limited by size. Using the graphical intuition of
Figure \ref{fig:rdd:var}, we provide a simpler proof of the same facts, with
the weak distance instead of the TV metric. Moreover, we add that
confidence intervals produced from Wald tests have zero confidence level.
It is worth highlighting the work of statisticians \cite{low1997} and \cite{cai2004}
on the impossibility of adaptation gains 
for confidence intervals on linear functionals of non-parametric functions. 
These authors take confidence intervals with correct coverage over a class of models 
$\mmp$ and derive a lower bound for the expected length of any confidence interval
under a given model $P \in \mmp$.
As the sample size increases, the rate at which these bounds shrink to zero does not depend on $P$. In other words, any confidence interval whose expected length at $P \in \mmp$ shrinks to zero faster than the lower bound 
must have incorrect coverage over $\mmp$.
 \cite{armstrongkolesar2015} 
derive a lower bound for the expected length of any confidence set that has correct coverage over  
$\mmp$. 
The lower bound increases to infinity as $\mmp$ becomes more general,
which is our impossibility type B.

On a positive note, the two types of
impossibility vanish if we restrict the class of models $\mmp$.
The approximation used to prove Corollary \ref{coro:rdd} fails if we assume that functions 
in $\m{G}$ have absolute slopes bounded by a finite constant $C$ on either side of the cutoff.
\cite{kamat2015} shows that Wald tests have correct size asymptotically if 
the first three derivatives of $g(x)$, as well as conditional moments, 
are uniformly bounded across $\mmp$. 
\cite{armstrongkolesar2015} derive minimax optimal-length confidence intervals
for a convex function class $\m{G}$  covering most smoothness or
shape assumptions used in econometrics.
In the RDD case, they consider functions $g(x)$ such that the $p^{th}$-order Taylor approximation residual  is bounded by $Cx^p$ on either side of the cutoff. 
In summary, applied researchers should bear in mind that the validity of tests and confidence sets
for the value of the discontinuity at the threshold relies heavily
on restricting the variation of average outcomes wrt the forcing variable $X$.
For example, consider the analysis  of summer-school programs in Chicago by \cite{jacob2004}.
The forcing variable $X$ is a standardized reading score determining eligibility for the program,
and $Y$ is a standardized test score in math or reading after the program.
Looking at their Figures 6 and 7, it seems reasonable to assume that the slope of the conditional mean of $Y$ given $X$ is smaller than one. In other words, an increase in today's reading score by 1 point increases tomorrow's average scores by less than one point.

Restricting the class of models $\mmp$ is not the only way to construct valid tests in RDD.
Another way to approach the problem is to consider null sets $\mmp_0$ different than those in 
Corollary \ref{coro:rdd}, where the focus is on the jump discontinuity at the threshold.
One example is the null hypothesis that an individual's outcome is solely affected by the treatment he receives and not by the effect of the treatment on neighboring individuals. 
In the summer-school application, the number of students attending classes in the summer is much smaller than during the school year.
It is likely that students in the summer program interact much more with each other, which leads to spillover effects of the treatment.
A researcher who desires to test for no spillovers specifies the null hypothesis of independence of $Y_i$ and $Y_j$
conditional on $X$. 
Under uniformly bounded moment conditions, \cite{shah2019} propose a test for such a null hypothesis that controls size and has non-trivial power in large samples.    

Another example of a null hypothesis that is immune to testing impossibility is
when absence of treatment effects is equivalent to  a smooth conditional mean function.
We may define the null hypothesis that $g$ is Lipschitz continuous with some constant $C$,
and the alternative hypothesis that $g$ is any other function as in Equation (\ref{eq:app:rdd:familyp}). 
Settings like this arise when the treatment variable $D$ is a function of the forcing variable $X$,
and this function changes at a known cutoff. 
This is the case of unemployment benefits in Austria, studied by \cite{card2015} (CLPW15).
For unemployed individuals that formerly earned $X$ less than a threshold $c$, the unemployment benefit grows with their earnings;
 otherwise, if they formerly earned  more than $c$, they simply receive a fixed benefit regardless of their earnings.
CLPW15 find that the unemployment duration does not depend on past earnings for those whose benefit is fixed to the right of the cutoff (see their Figure 3).
Moral hazard leads to unemployment duration that increases as benefits increase with income to the left of the cutoff.
Therefore, the researcher may specify the null hypothesis of a smooth conditional mean to test for the lack of  moral hazard. 
Rejections may occur because of a sudden change in slope or a jump discontinuity at the threshold, both of which are evidence of a change in behavior regarding job search.
Note, however, that the null hypothesis of Lipschitz $g$
is different than the null hypothesis in the so-called Regression
Kink Design (RKD) studied by CLPW15. 
The RKD null states that the first derivative of $g$ is continuous at the threshold, 
and this null suffers from testing impossibility.

RKD has recently gained popularity in economics. In addition to CLPW15, see  
\cite{dong2014}, \cite{nielsen2010}, and \cite{simonsen2016}.
 The
setup is the same as in the RDD case, except that the causal effect of
interest is the change in the slope of the conditional mean of outcomes at
the threshold. Continuity of the first derivatives $\nabla _{x}\mathbb{E}%
[Y_{i}(1)|X_{i}=x]$ and $\nabla _{x}\mathbb{E}[Y_{i}(0)|X_{i}=x]$ at the
threshold $x=c$ guarantees identification of the average effect. The
parameter of interest $m=\mu (P)$ is a function of the distribution of $%
Z_{i}=(X_i, Y_{i})$: 
\begin{equation}
\mu (P)=\nabla _{x}\mathbb{E}[Y_{i}(1)-Y_{i}(0)|X_{i}=x]=\lim\limits_{x%
\downarrow c}\nabla _{x}\mathbb{E}[Y_{i}|X_{i}=x]-\lim\limits_{x\uparrow
c}\nabla _{x}\mathbb{E}[Y_{i}|X_{i}=x].
\end{equation}

The family of all possible distributions of $Z_i$ is defined in a slightly
different way than in Equation (\ref{eq:app:rdd:familyp}): 
\begin{equation}
\mathbb{P} = \{P: (X_i, Y_i) \sim P,~ \exists g \in \mathcal{G} \text{ s.t. }
\nabla_x \mathbb{E}[Y_i|X_i=x] = g(x) \}.
\end{equation}

Weak identification of $\mu$ arises from the fact that any conditional mean
function $\mathbb{E}[Y_i|X_i=x]$ with a discontinuous first derivative at $%
x=c$ is well-approximated by a sequence of continuously differentiable
conditional mean functions. Assumption \ref{assu:dense} is easily verified
using this insight.

\begin{corollary}
\label{coro:rkd} Assumption \ref{assu:dense} is satisfied for $\mmp_{0,m}$ $\forall m\in 
\mathbb{R}$, and Theorems \ref{theo:impossible_test} and \ref{theo:ci} apply
to RKD. Namely, (i) a.s. continuous tests $\phi _{m}(Z)$ for the value of the
kink discontinuity $m$ have power limited by size; and (ii) confidence sets for
the value of the kink discontinuity $m$ and with finite expected length have
zero confidence level.
\end{corollary}

The proof of Corollary \ref{coro:rkd} follows that of Corollary \ref{coro:rdd}.
 Simply use the new definitions of $\mathbb{P}$
and $\mu (P)$, and construct the sequence $P_{k}$ with $\nabla _{x}\mathbb{E}_{P_{k}}[Y_{i}|X_{i}=x]=g_{k}(x)$.

\subsection{Testing for the Existence of Bunching}
\label{sec:app:bun}

\indent

The second example applies Theorem \ref{theo:impossible_test} to the problem of testing for the existence of bunching in a scalar random variable. 
Bunching occurs when the distribution of $X$ exhibits a non-zero probability at known point $x_0$, but it is continuous in a neighborhood of $x_0$.
Bunching in the distribution of a single variable is the object of interest in many empirical studies.
For example, \cite{saez2010} and \cite{kleven2013} rely on the existence of bunching on ``reported income'' at the boundary of tax brackets to identify the elasticity of reported income wrt tax rates;
\cite{goncalves2017} use bunching on ``charged speed in traffic tickets'' to separate lenient from non-lenient police officers and identify racial discrimination;
and a standard practice in RDD analyses is to check if the distribution of the forcing variable has bunching at the cutoff, which would count as evidence against the design.

Suppose $X$ is a scalar random variable. 
In the absence of bunching, assume the CDF of $X$ is continuously differentiable.
Testing for bunching amounts to testing whether $X$ has positive probability mass at $x_0$.
Let $\mmp_0$ be the set of distributions of $X$ with a continuously differentiable CDF. 
The set $\mmp_1$ is all mixed continuous-discrete distributions, with one mass point at $x_0$, but continuously differentiable CDF otherwise.\footnote{The
assumption that the CDF is continuously differentiable is not necessary in this section.
We impose this assumption because typical non-parametric density estimators assume a continuous density.
The testing impossibility of this section occurs regardless of whether the CDF is assumed continuously differentiable, or simply continuous.}
Any distribution $Q$ under the alternative is well-approximated in the LP metric by a sequence of distributions $P_k$ under the null.
Therefore, any a.s. continuous test has power limited by size.
\begin{corollary}
\label{coro:bun}
Assumption \ref{assu:dense} is satisfied in the problem of testing for the existence of bunching.
Hence, any test $\phi(Z)$ that is a.s. continuous under $\mmp_1$ has power limited by size.
\end{corollary}

There is one interesting feature about this example that is not shared by the RDD and RKD examples of the previous section.
In this example, it is not possible to find a sequence $P_k$ under the null that approximates a  $Q \in \mmp_1$ using the TV metric. 
The event $X=x_0$ always has zero probability under the null, but strictly positive probability under the alternative.
Therefore, $d_{TV}(P,Q) >0$ for every $P \in \mmp_0, Q \in \mmp_1$. Theorem 1 suggests that there exists a test
whose maximum power is bigger than size, but our Theorem \ref{theo:impossible_test} says this test cannot be a.s. continuous under $\mmp_1$.

The use of the LP metric, as opposed to the TV metric, leads us to search for tests that are discontinuous under $\mmp_1$.
For a sample with $n$ iid observations $X_i$, the test $\phi(X_1,\ldots,X_n) = \mmi \left\{ \frac{1}{n}\sum_{i=1}^n \mmi\{X_i = x_0\} >0 \right\}$
is discontinuous under $\mmp_1$. This test has size equal to zero, and power equal to $1 - (1-\delta)^n$,
 where $\delta = \mmp[X_i=x_0]$.

\subsection{Exogeneity Tests Based on Bunching}

\indent

The third example comes from \cite{caetano2015}, who uses the idea of bunching in a conditional distribution of $Y$ given $X$ to construct an exogeneity test that does not require instrumental variables.
It applies to regression models where the distribution of unobserved  factors is assumed to be discontinuous wrt an explanatory variable.
Of interest is the impact of a scalar explanatory variable $X$ on an outcome variable $Y$, after controlling for covariates $W$.
For example, suppose we are interested in the effect of average number of cigarettes smoked per day $X$ on birth weight $Y$, after controlling for mothers' observed characteristics $W$. 
Conditional on $(X,W)$, the distribution of mothers' unobserved characteristics $U$ is said to bunch at zero smoking if it changes drastically when we compare non-smoking mothers to mothers that smoke very little. 
If bunching occurs, then the variable $X$ is endogenous because we cannot separate the effect of smoking on birth weight from the effect of unobserved characteristics on birth  weight.

The population model that determines $Y$ is written as $Y=h(X,W)+U$, where $%
U $ summarizes unobserved confounding factors affecting $Y$. 
We are unable to infer bunching on $U$ unless $h$ is assumed continuous on $(X,W)$. 
Bunching of $U$ wrt $X$ is evidence of local endogeneity of $X$ at $X=0$. 
Bunching at 0 implies discontinuity of $\mathbb{E}[U|X=0,W]-\mathbb{E}[U|X=x,W]$ as $x \downarrow 0$.
Continuity of $h$ makes bunching equivalent to a discontinuity of 
$\mathbb{E}[Y|X=0,W=w]-\mathbb{E}[Y|X=x,W=w]$ as $x \downarrow 0$ for every $w$. \cite{caetano2015}
proposes testing 
\begin{equation}
\forall w~\lim_{x \downarrow 0 } \mathbb{E}[Y|X=0,W=w]-\mathbb{E}[Y|X=x,W=w]=0
\end{equation}%
as a means of testing for local exogeneity of $X$ at $X=0$. We argue that $h$
may have a high slope on $X$, or even be discontinuous on $X$, which makes
exogeneity untestable.

The observed data $Z=(Z_{1},\ldots ,Z_{n})$, $Z_{i}=(X_{i},W_{i},Y_{i})$ is
iid with probability $P$. The support of $(X_{i},W_{i})$ is denoted $\mathbb{X}\times 
\mathbb{W}$.
The distribution of $Y$ conditional on $(X,W)$ is assumed to be continuous.
The distribution of $X$ has non-zero probability at $X=0$, but it is continuous otherwise.
Assume $\exists \delta >0$ such that $[0,\delta )\subset 
\mathbb{X}$. 

Let $\mathcal{G}$ denote the space of all functions $g:\mathbb{X}\times \mathbb{W}\rightarrow \mathbb{R}$ that are bounded
and infinitely many times continuously differentiable wrt $x$ over $\{ \mathbb{X}\setminus \{0\}\} \times \mmw$. 
The size of the discontinuity at $X=0$ may take any value in $\mmr$.
The family of all
possible distributions is denoted as 
\begin{equation}
\mathbb{P}=\{P:Z_{i}\sim P,~\exists g\in \mathcal{G}\text{ s.t. }\mathbb{E}%
_{P}[Y_{i}|X_{i}=x,W_{i}=w]=g(x,w)\}.
\end{equation}

Under local exogeneity of $X$, the function $\tau _{P}(w)=\mathbb{E}%
_{P}[Y_{i}|X_{i}=0,W_{i}=w]-\lim\limits_{x\downarrow 0}\mathbb{E}_P[Y_{i}|X_{i}=x,W_{i}=w]$ must be equal to $0~\forall w\in \mathbb{W}$. In
practice, it is convenient to conduct inference on an aggregate of $\tau
_{P}(w)$ over $w\in \mathbb{W}$ instead of on the entire function $\tau
_{P}(w)$. Examples of aggregation include the average of $|\tau _{P}(W)|$,
the square root of the average of $\tau _{P}(W)^{2}$, or the supremum of $|\tau _{P}(w)|$ over $w \in \mathbb{W}$.
For the sake of brevity, we choose the second option. 
For a distribution $P\in \mathbb{P}$, define $\mu (P) = \left[ \mme_P\left( \tau _{P}(W)^2 \right) \right]^{1/2}$.
Local exogeneity corresponds to the test of $\mu (P)=0$ versus $\mu (P)\neq 0$.

The parameter $\mu (P)$ is weakly identified in the class of models $\mathbb{%
P}$. Just as in the RDD case, any conditional mean function $\mathbb{E}%
[Y_{i}|X_{i}=x,W_{i}=w]$ with a discontinuity at $x=0$ is well-approximated
by a sequence of continuous conditional mean functions $\mathbb{E}%
[Y_{i}|X_{i}=x,W_{i}=w]$. Assumption \ref{assu:dense} is verified using the
same argument as in the RDD case.

\begin{corollary}
\label{coro:exo} Assumption \ref{assu:dense} is satisfied for $\mmp_{0,m}$ $\forall m\in 
\mathbb{R}$, and Theorems \ref{theo:impossible_test} and \ref{theo:ci} apply
to the case of the local exogeneity test. Namely, (i) a.s. continuous tests $%
\phi _{m}(Z)$ for the value of the aggregate discontinuity $m$ have power
limited by size; and (ii) confidence sets for the value of the aggregate
discontinuity $m$ and with finite expected length have zero confidence level.
\end{corollary}

The inference procedures suggested by \cite{caetano2015} rely on
non-parametric local polynomial estimation methods. As in the RDD case, these
procedures yield tests that are a.s. continuous in the data and confidence
intervals of finite expected length. 
Corollary \ref{coro:exo} implies lack of size control and zero confidence level.


\subsection{Time-Series Models}
\label{sec:app:ts}

\indent

The fourth example illustrates robust hypothesis testing wrt the LP metric, and it is of practical relevance to macroeconomists.
Macroeconometrics often uses linear time-series processes.
This is motivated by the Wold  Representation Theorem, which asserts that every covariance-stationary process $x_{t}$ can be written as an MA process plus some
deterministic term:
\begin{equation*}
x_{t}=B\left( L\right) \eps _{t},
\end{equation*}%
where $L$ is the lag operator, $B\left( l\right) =\sum_{i=0}^{q}b_{i}~l^{i}$, and $\eps _{t}$ is an uncorrelated error sequence.
A caveat is that the order $q$ needs to be too large to be useful for many applications.
The features of MA processes with infinite lag order are well captured by ARMA models
\begin{equation*}
A\left( L\right) x_{t}=B\left( L\right) \eps _{t},
\end{equation*}%
with small orders for $A\left( L\right) $ and $B\left( L\right) $, where $A\left( l\right) =\sum_{i=0}^{p}a_{i} ~l^{i}$.
The closure of the set of stationary ARMA(p,q) models with finite order $(p,q)$ does not necessarily contain only stationary models.
The simplest example happens when $A\left( l \right) =1 - a  l$ and $B\left( l\right) =1$. 
The process is stationary when $|a|<1$, but it is non-stationary when $a=1$.
 This observation led to ARIMA models, which better capture the persistence in time series.

Starting in the 1990s, applied researchers began to realize that ARIMA
models themselves have limitations. This led to the development of other
stochastic processes, including error duration models, Markov switching
models, threshold models, structural breaks, and fractionally integrated
processes, among others. This is a vast literature and includes papers by
\cite{hamilton1989}, \cite{parke1999}, and \cite{bai1998}, just to name a
few.

A number of authors point out that these different model extensions may not be too far from each other.
For example, \cite{perron1989} shows that integrated processes with drift and stationary models with a broken trend can be easily confused;
\cite{parke1999} points out that the error-duration model encapsulates fractionally integrated series;
\cite{granger1999} and \cite{diebold2001} find that linear processes with breaks can be misinterpreted as long-memory models.
In these papers, and in most of the related econometrics literature, the focus is on the autocovariance of the stochastic process.

Our discussion of robust hypothesis testing in Section \ref{sec:rob} suggests looking at the closure of ARMA processes to distinguish these models from each other.
For example, take the problem of testing the null that a process is covariance-stationary, against the alternative that it is an error-duration model or a Compound Poisson model.
The existence of a test with non-trivial power requires us to look for the TV distance between these sets of processes. 
However, the ability to approximate theses processes in the TV distance is often based on quite stringent assumptions.
For example, see \cite{barbour1999} for the TV approximation of Compound Poisson processes.

The problem of searching for tests for covariance-stationary versus error-duration or Compound Poisson becomes much easier if we focus on the LP metric.
To solve this problem, we rely on \cite{bickel1996}, whose work has been largely ignored in the econometrics literature.
They characterize the closure of AR and MA processes wrt the TV and the Mallows metric (also known as the Wasserstein metric).
The TV metric is stronger than the Mallows metric, which in turn is stronger than the LP metric.
Indeed, convergence under the Mallows metric implies weak convergence and convergence in second moments; see \cite{bickel1981} and \cite{bickel1996}.
As a result, the closure of stochastic processes wrt the LP metric is larger than the closure wrt the Mallows metric.
It turns out that error-duration and Compound Poisson models are in the closure wrt the LP metric.
In other words, the robustified null set wrt the LP metric contains the alternative set, and the minimal TV distance between these sets is zero.
Hence, all tests for the robustified null have power no larger than size.

Given that the closure of ARMA processes of infinite order is quite rich, we may wonder which hypotheses are testable.
\cite{bahadur} and \cite{romano2004} point out  it is hopeless to test population means, even in the iid case without further moment constraints.
Could we try to test quantiles? 
\cite{peskir2000} and \cite{shorack2009}  provide sufficient conditions for uniform convergence of empirical processes under time dependence.
A natural choice for quantile testing is the value at risk (VaR), which is commonly used in the finance literature.
It would be interesting to establish the class of empirical processes for which hypotheses for the VaR are testable.
We leave this example for future work.

\section{Simulations}
\label{sec:simul}

\indent 

In this section, we provide Monte Carlo simulations to illustrate the impossibility of testing within the context of RDD. 
We find that the Wald test fails to control size uniformly under the null hypothesis.
We use a data-generating process (DGP) based on an empirical example. 
Lack of size control occurs even for DGPs that are consistent with the data.
Moreover, the simulations also show that the Wald test has very little power after artificially controlling  size.
For the sake of brevity, we focus on the RDD case, and we expect similar findings for the RKD and Exogeneity Test cases.

\begin{figure}[H]
\caption{Conditional Mean Function Based on Lee's (2008) Data }
\label{fig:rdd:lee}
\vspace{-2 \baselineskip}  
\begin{center}
\includegraphics[width=5in]{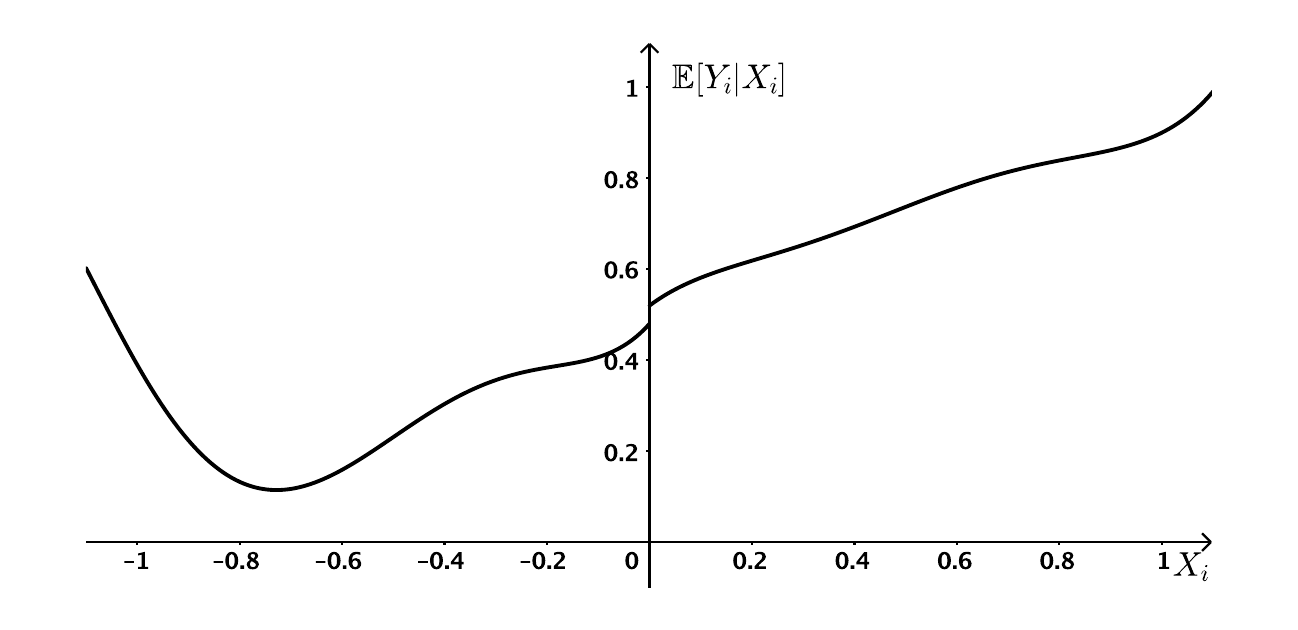}
\end{center}
\vspace{-.9cm}
\caption*{\footnotesize Notes: conditional mean function of Equation (\ref{dgp:rdd:1}). The forcing variable $X$ is the margin of victory of the Democratic party in percentage of votes in time $t$. The outcome variable $Y$ is equal to one if Democrats win in time $t+1$, but equal to zero otherwise.   }
\end{figure}

Our DGP is based on the incumbency data of \cite{lee2008}.
Lee studies incumbency advantage in the US House of Representatives.
Districts where a party's candidate barely wins an election are, on average, comparable to districts where that party's candidate barely loses the election. 
The forcing variable $X$ is the margin of victory of the Democratic party in percentage of votes. 
The target parameter is the effect of the Democrats winning  the election at time $t$ (incumbency) on the probability of the Democrats winning the election at time $t+1$.
Lee's data have been used for simulation studies by several other econometricians, for example, by \cite{imbens2012optimal}, \cite{cattaneo2014calonico}, and \cite{armstrongkolesar2015}.
We use the Monte Carlo DGP of \cite{imbens2012optimal} and \cite{cattaneo2014calonico},
described in Equation (\ref{dgp:rdd:1}).\footnote{The 
DGP in Equation (\ref{dgp:rdd:1}) belongs to the class of functions that \cite{armstrongkolesar2015}
study in their application to RDD. 
The set of functions $\m{F}_{RDP,p}(C)$ on their page 658 contains 
Equation (\ref{dgp:rdd:1}) with $p=2$ and constant $C=7.2$.}
\begin{gather}
Y=
\left\{
\begin{array}{ll}
 0.48  + 1.27 X + 7.18 X^2 + 20.21 X^3 & 
 \\
 \hspace{4cm} + 21.54 X^4 + 7.33 X^5  + U& \text{ if } X\in (-0.99,0)
 \\
 \\
 0.52  + 0.84 X - 3 X^2 + 7.99 X^3 & 
 \\
 \hspace{4cm} - 9.01 X^4 + 3.56 X^5 +U & \text{ if } X \in [0, 0.99]
 \end{array}
 \right.
 \label{dgp:rdd:1}
\end{gather}
where  $X$ is distributed as Beta$(2,4)$, $U$ is zero-mean Gaussian with standard deviation $0.1295$, and $X$ is independent of $U$. Figure \ref{fig:rdd:lee} depicts the conditional mean function of Equation (\ref{dgp:rdd:1}).
\vspace{-1 \baselineskip}  

Our simulation study uses variations of  Equation (\ref{dgp:rdd:1}) that are governed by two parameters: $\tau \in \mmr$ and $M \in \mmr_+$.
\begin{gather}
Y=
\left\{
\begin{array}{ll}
 0.48  + \tau \Lambda\left( 4 M X / \tau \right) + 1.27 X + 7.18 X^2 + 20.21 X^3 & 
 \\
 \hspace{4cm} + 21.54 X^4 + 7.33 X^5  + U& \text{ if } X\in (-0.99,0)
 \\
 \\
 0.48  + \tau \Lambda\left( 4 M X / \tau \right) + 0.84 X - 3 X^2 + 7.99 X^3 & 
 \\
 \hspace{4cm} - 9.01 X^4 + 3.56 X^5 +U & \text{ if } X\in [0,0.99)
 \end{array}
 \right.
 \label{dgp:rdd:2}
\end{gather}
where  $\Lambda\left( \cdot \right)$
is the logistic CDF function.

\begin{figure}
\caption{Approximating a Discontinuous Conditional Mean Function ($\tau=0.04$)}
\label{fig:rdd:var}
\vspace{-.75 \baselineskip}  
\begin{center}
  \begin{minipage}{.5 \textwidth }
    \centering
    (a) $M=0$
    
    \includegraphics[width=3.0in]{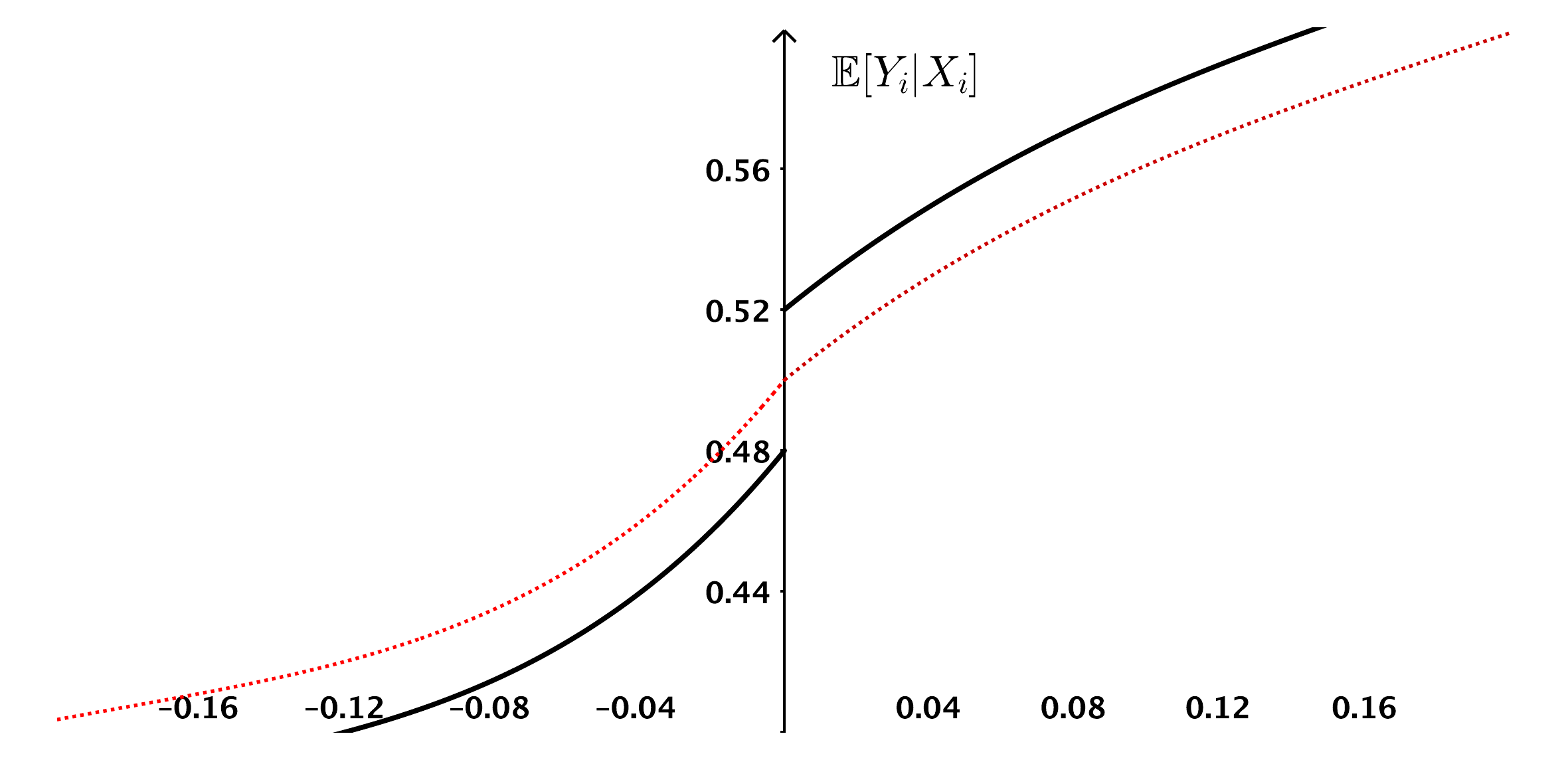}
  \end{minipage}%
  \begin{minipage}{.5 \textwidth }
    \centering
    (b) $M=0.5$

    \includegraphics[width=3.0in]{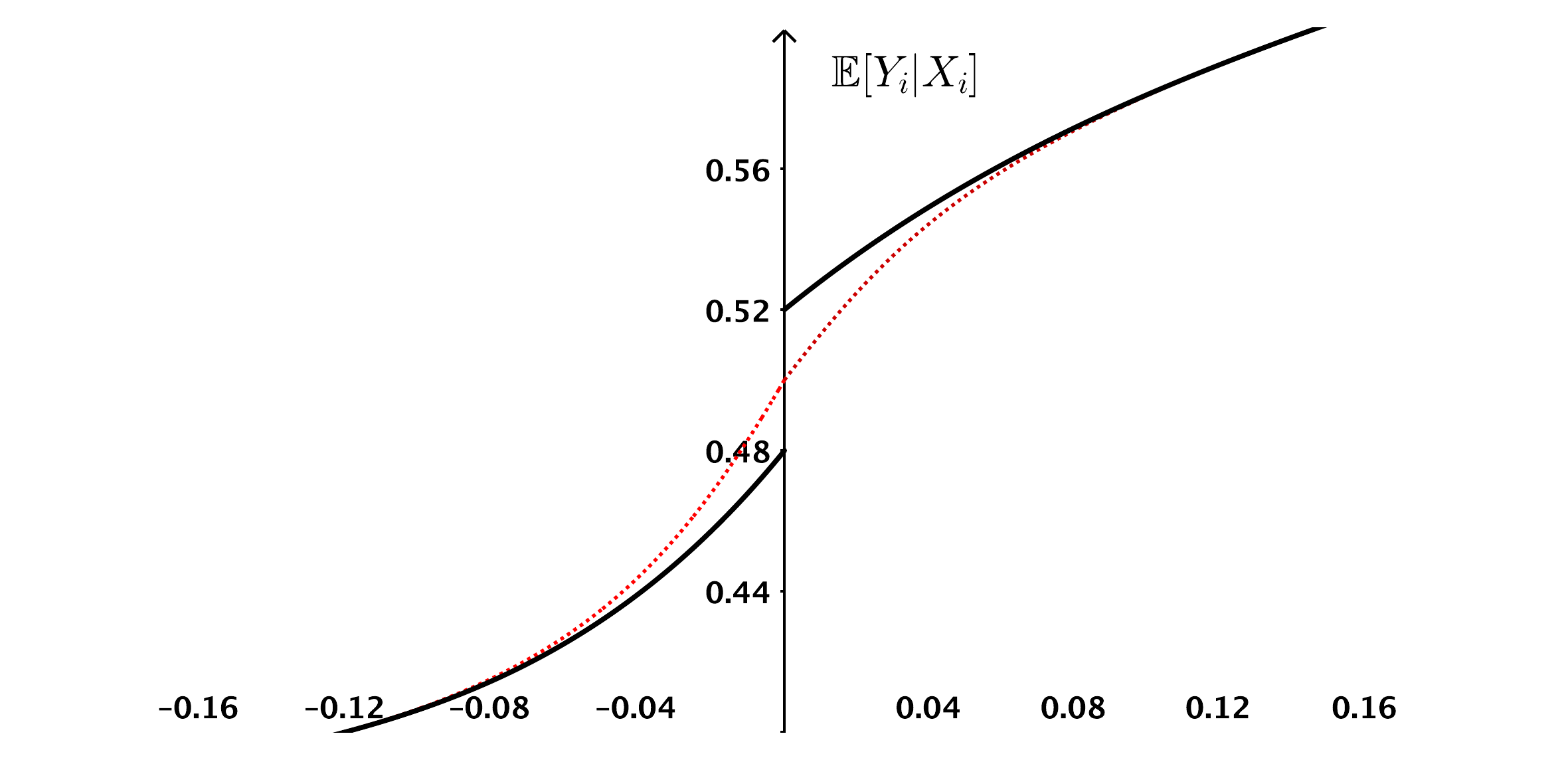}
  \end{minipage}%
  
  \bigskip
  
  \begin{minipage}{.5 \textwidth }
    \centering
    (c) $M=2$

    \includegraphics[width=3.0in]{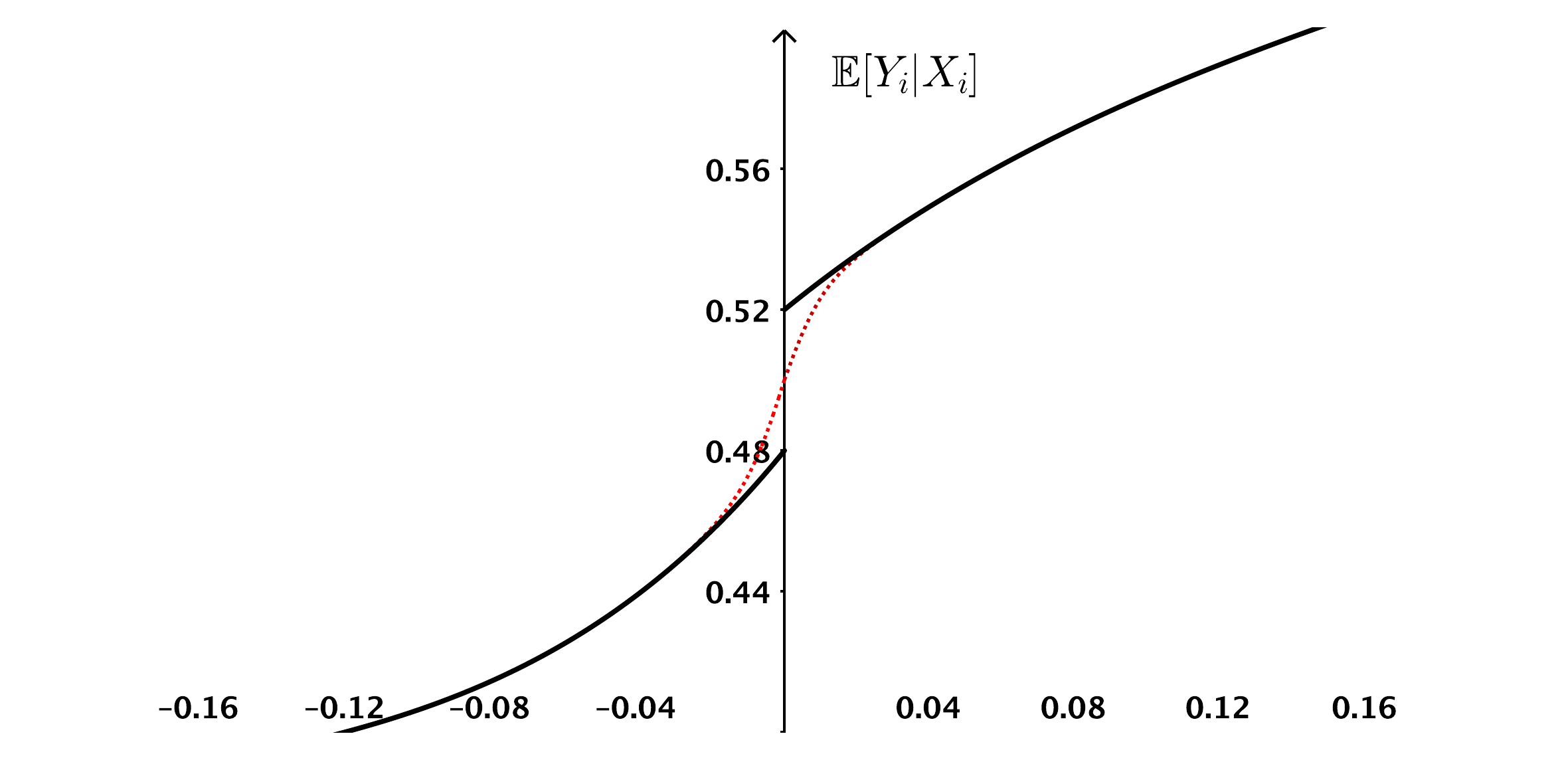}
  \end{minipage}%
  \begin{minipage}{.5 \textwidth }
    \centering
    (c) $M=8$

    \includegraphics[width=3.0in]{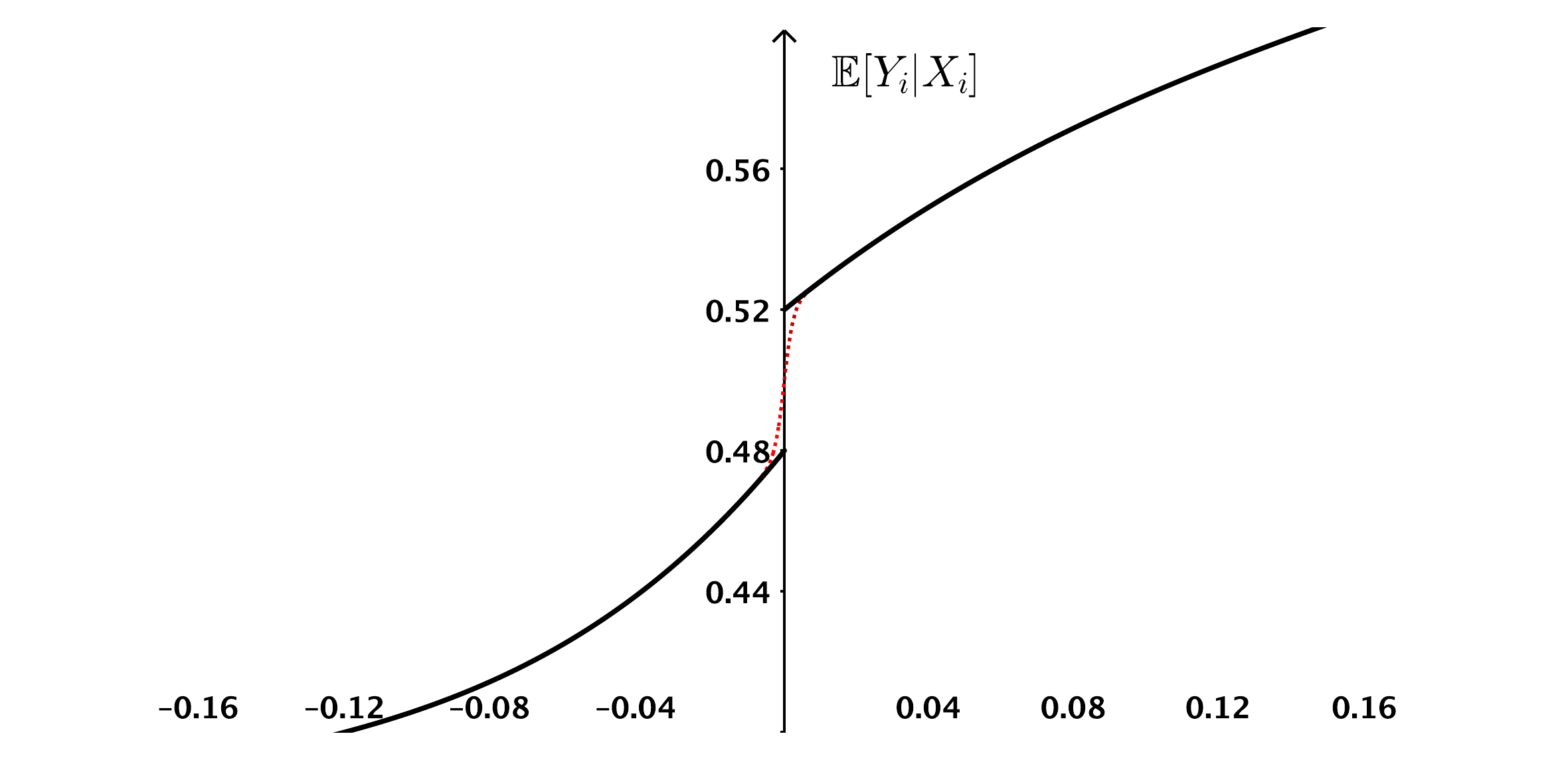}
  \end{minipage}%

\end{center}
\vspace{-.5cm}
\caption*{\footnotesize
Notes: 
the discontinuous conditional mean function $\mme[Y|X]$ (solid line) is approximated by a sequence of continuous conditional mean functions (dotted lines).
The solid line is the $\mme[Y|X]$ of Model \ref{dgp:rdd:1},
and the dotted line is the $\mme[Y|X]$  of Model \ref{dgp:rdd:2}
for $\tau=0.04$ and $M\in\{ 0, 0.5, 2, 8\}$. 
The figure illustrates that 
model \ref{dgp:rdd:2} approximates the DGP based on \cite{lee2008}
as the slope at $X=0$ grows large.
}
\end{figure}

The conditional mean function of both Equations (\ref{dgp:rdd:1}) and (\ref{dgp:rdd:2})
are differentiable on either side of the cutoff.
The first is discontinuous at $X=0$ with discontinuity of size $0.04$, 
while the second is continuous at $X=0$.
For $\tau=0.04$,
Equation (\ref{dgp:rdd:2}) approximates Equation (\ref{dgp:rdd:1}) as $M \to \infty$.
The parameter $M$ is the derivative of $ \tau \Lambda\left( 4 M X / \tau \right)$ wrt $X$ at $X=0$.
As the slope $M$ grows large, the continuous conditional mean function of Equation (\ref{dgp:rdd:2}) approximates a discontinuous function with discontinuity of size $\tau$.
 Figure \ref{fig:rdd:var} illustrates this approximation, as well as the proof of Corollary \ref{coro:rdd} in Section \ref{sec:app:rd}. 
For example, a model similar to Equation (\ref{dgp:rdd:2}) with high values of $M$ arises when districts manipulate the share of votes in order to win the election.
Manipulation of the forcing variable has been extensively studied in the RDD literature. 
See, for example, \cite{mccrary2008} and \cite{gerard2016}.
Suppose the average causal effect of winning the election conditional on $X$ is small for districts with small margin of victory,
but large otherwise. 
In the absence of manipulation, $\mme[Y|X]$ is continuous and very smooth to the right of the cutoff.
The party in  districts with low $X$ has incentives to manipulate the election, and the researcher observes 
the manipulated margin of victory $\tilde{X}$, instead of $X$. 
Assume the probability that manipulation occurs conditional on $\tilde{X}$ increases continuously but sharply to the right of the cutoff.
In this case, the researcher observes a conditional mean function  that is continuous at the cutoff but that increases sharply to the right of the cutoff.
In practice, one may falsely reject the null of zero effect simply because of manipulation, and not because of an actual causal effect.
We provide a concrete example for this DGP in Section \ref{supp:rdd:manipul} in the appendix.

The parameter of interest is $m$, the size of the jump discontinuity at $X=0$. 
The null hypothesis is $m=0$, which is the set of models in Equation (\ref{dgp:rdd:2}) with $\tau \in \mmr$ and  $M \in \mmr_+$.
The alternative hypothesis is $m \neq 0$, which is the set of models with  $\tau \neq 0$ and  $M = \infty$.
Section \ref{sec:app:rd} shows that any model in the alternative is well-approximated in the LP metric by models under the null. 
The power of a.s. continuous tests is less than or equal to size.

The Monte Carlo experiment simulates 10,000 draws of an iid sample with 500 observations.
The range of $(\tau,M)$ values for Model \ref{dgp:rdd:2} in the experiment is consistent with the magnitudes of Lee's DGP.
The maximum slope magnitude of the conditional mean graph in Figure \ref{fig:rdd:lee} is 1.97,
and we set $M \in \{0, 2,  \ldots, 10 \}$.
The value of $m$ for Lee's DGP is $0.04$, and we vary $\tau$ in $\{0, 0.01, 0.04, 0.08\}$.
We conduct a size and a power analysis.
In the size analysis, we simulate rejection probabilities of the Wald test
 under each  $(\tau,M)$-model.
  The estimates of $m$ and standard errors are obtained by the robust bias-corrected method of 
\cite{cattaneo2014calonico} and implemented using the STATA package \verb+rdrobust+.
For each model $(\tau,M)$, the critical value of the test comes from the simulated distribution of the statistic under model $(\tau,0)$. 
This ensures exact size of the test in the smoothest model under the null $(M=0)$.
\begin{table}
\begin{center}
\caption{Rejection Probability Under the Null - Size 5\%}
\label{table:rdd:size:5}
 \begin{tabular}{c  c  c  c  c  c  c  } 
 \hline \hline 
$\tau$  & $ M= 0 $  & $ M= 2 $  & $ M= 4 $  & $ M= 6 $  & $ M= 8 $  & $ M= 10 $ 
\\
\hline
 .01  &   0.0500 &   0.0540 &   0.0557 &   0.0592 &   0.0585 &   0.0580 \\ 
 .02  &   0.0500 &   0.0665 &   0.0678 &   0.0649 &   0.0694 &   0.0685 \\ 
 .03  &   0.0500 &   0.0910 &   0.0942 &   0.0938 &   0.0941 &   0.1016 \\ 
 .04  &   0.0500 &   0.1005 &   0.1067 &   0.1071 &   0.1121 &   0.1139 \\ 
 .05  &   0.0500 &   0.1114 &   0.1264 &   0.1334 &   0.1464 &   0.1434 \\ 
 .06  &   0.0500 &   0.1292 &   0.1632 &   0.1680 &   0.1819 &   0.1819 \\ 
 .07  &   0.0500 &   0.1258 &   0.1617 &   0.1832 &   0.1906 &   0.2026 \\ 
 .08  &   0.0500 &   0.1320 &   0.1888 &   0.2142 &   0.2266 &   0.2427 \\ 
\hline  \hline 
\end{tabular} 

\caption*{\footnotesize
Notes: 
the table displays the simulated rejection probability of the Wald test under various choices of $(\tau,M)$ for Model \ref{dgp:rdd:2}.
Critical values of the test vary by row, but are constant across columns. 
For each $(\tau,M)$-model, the critical value of the test comes from the simulated distribution of the statistic under model $(\tau,0)$.
The estimates of $m$ and standard errors for the Wald test are obtained by the robust bias-corrected method of 
\cite{cattaneo2014calonico} and implemented using the STATA package ``rdrobust''.
}
\end{center}
\end{table}

The nominal size of the Wald tests in Table \ref{table:rdd:size:5} is 5\%, and the simulated rejection probability increases with $\tau$ and $M$.
For the maximum slope of $M=2$ observed from the model in Equation (\ref{dgp:rdd:1}),
the size of the test varies between 5.4\% and 13.2\%, depending on the choice of the model under the null. The true value of $M$ is unknown, and a more conservative upper bound on the slope $M=10$ distorts the size of the test up to 24\%.

In the power analysis, we study rejection probabilities for models with $M=\infty$ and $\tau \in \{0, 0.01, \ldots, 0.08\}$. 
These models fall under the alternative because $m=\tau$ when $M=\infty$.  
For each $(\tau,\infty)$-model, we would like the test to have correct size under the least favorable null model.
Table \ref{table:rdd:size:5} suggests that the least favorable model under the null is the one with the highest slope $M$. 
Figure \ref{fig:rdd:var} shows that null models can approximate any alternative $(\tau,\infty)$-model arbitrarily well.
If we restrict the slope at $X=0$ to be at most $M$, the worst-case model under the null for the alternative  $(\tau,\infty)$-model
is the  $(\tau,M)$-model. 
To evaluate the rejection probability under a $(\tau,\infty)$-model, the critical value of the test comes from the simulated distribution of the statistic under a $(\tau,M)$-model  for various choices of $(\tau,M)$.
That way, the test has correct size when $m=0$ under all possibilities of least favorable $(\tau,M)$-models. 
\begin{table}[H]
\begin{center}
\caption{Rejection Probability Under the Alternative - Size 5\%}
\label{table:rdd:power:5}
 \begin{tabular}{c  c  c  c  c  c  c  } 
 \hline \hline 
$\tau$  & $ M= 0 $  & $ M= 2 $  & $ M= 4 $  & $ M= 6 $  & $ M= 8 $  & $ M= 10 $ 
\\
\hline
 .01  &   0.0610 &   0.0508 &   0.0504 &   0.0501 &   0.0504 &   0.0500 \\ 
 .02  &   0.0763 &   0.0527 &   0.0524 &   0.0505 &   0.0513 &   0.0501 \\ 
 .03  &   0.1020 &   0.0574 &   0.0532 &   0.0525 &   0.0526 &   0.0527 \\ 
 .04  &   0.1204 &   0.0646 &   0.0571 &   0.0556 &   0.0536 &   0.0524 \\ 
 .05  &   0.1583 &   0.0770 &   0.0618 &   0.0597 &   0.0573 &   0.0544 \\ 
 .06  &   0.2013 &   0.0899 &   0.0682 &   0.0638 &   0.0605 &   0.0569 \\ 
 .07  &   0.2192 &   0.1023 &   0.0732 &   0.0677 &   0.0642 &   0.0590 \\ 
 .08  &   0.2781 &   0.1179 &   0.0839 &   0.0707 &   0.0654 &   0.0631 \\ 
\hline  \hline 
\end{tabular} 

\caption*{\footnotesize
Notes: 
the entries of the table display the simulated rejection probability of the Wald test under Model \ref{dgp:rdd:2} 
with various $\tau$ and $M=\infty$, so that the size of the discontinuity is $m=\tau$.
Critical values of the test vary by row and column. 
For each $(\tau,M)$-entry, the critical value comes from the simulated distribution of the statistic under a null $(\tau,M)$-model.
The estimates of $m$ and standard errors for the Wald test are obtained by the robust bias-corrected method of 
\cite{cattaneo2014calonico} and implemented using the STATA package ``rdrobust''.
}
\end{center}
\end{table}
\vspace{-2 \baselineskip}  
  
The power of the tests in Table \ref{table:rdd:power:5} increases with the size of discontinuity $\tau$, but it decreases with the slope $M$ of the least favorable model under the null. 
Intuitively, the higher $M$ is, the harder it becomes to distinguish a $(\tau,M)$-model from a $(\tau,\infty)$-model. 
For the empirically relevant values of $\tau=0.04$ and $M=2$, we see that the power of the test is 6.5\%,
barely above its size.
More conservative upper bounds on the slope of the model under the null essentially make power equal size.
Section \ref{supp:simul:rdd} in the appendix contains versions of these tables for nominal levels 1\% and 10\%, as well as the simulated critical values used.

\section{Conclusion}

\label{sec:con}

\indent

When drawing inference on a parameter in econometric models, some authors
provide conditions under which tests have trivial power (impossibility type A).
 Others examine when confidence regions have error
probability equal to one (impossibility type B). 
The motivation
behind these negative results is that the parameter of interest may be
nearly unidentified across models. 
Impossible inference relies on models being indistinguishable wrt some notion of distance.
Some authors distinguish models using the Total Variation (TV) metric and others rely on the L\'{e}vy-Prokhorov (LP) metric, which is a weaker notion of distance.
The ability to distinguish models in the TV metric is a necessary and sufficient condition for the existence of tests with non-trivial power.
Impossible inference in terms of a weaker notion of distance is often easier to prove, it is applicable to the widely-used class of  almost surely continuous tests, and it is useful for robust hypothesis testing.

Impossibility type A is stronger than type B. 
\cite{dufour1997} focuses on models in which tests based on bounded confidence regions fail to control size, but they could still have non-trivial power.
Take the simultaneous equations model when instrumental variables may be arbitrarily weak. 
\cite{moreira2002,moreira2003} and \cite{kleibergen2005} propose tests that have correct size in models with type B impossibility.
Furthermore, these tests have good power when identification is strong, being efficient under the usual asymptotics.
Their power is not trivial, exactly because not every model under the alternative is approximated by models under the null.


The choice of the LP versus the TV metric connects our work to the work of Peter J. Huber on robust statistics.
It leads us to look at the closure of model departures under the LP metric.
In particular, robust hypothesis testing requires a non-zero TV distance between the closure of the null and alternative sets under the LP metric. 
For example, it is impossible to find a robust test that powerfully distinguishes 
covariance-stationary models from error-duration and Compound Poisson models, because the closure of the former contains the latter.
This closure is quite rich, and we wonder what sort of hypotheses are testable.
It is impossible to test the population mean, so one possibility may be quantiles such as value at risk (VaR).
\cite{peskir2000} and \cite{shorack2009}  provide sufficient conditions for convergence of empirical processes under dependence.
It would be interesting future work to build on these conditions to establish the class  of processes in which quantile testing is possible.

\section{Acknowledgements}

We thank Tim Armstrong, Leandro Gorno, and anonymous referees for helpful comments and suggestions.
Bertanha gratefully acknowledges support from ISLA-Notre Dame and CORE-UcLouvain.
Moreira
acknowledges the research support of CNPq and FAPERJ.
This study was financed in part by the Coordena\c{c}\~{a}o de Aperfei\c{c}oamento 
de Pessoal de N\'{i}vel Superior - Brasil (CAPES) - Finance Code 001.

\begin{singlespace}

\bibliographystyle{econ}
\bibliography{biblio}

\appendix

\numberwithin{lemma}{section}

\section{Appendix}

\subsection{Proof of Corollary \ref{coro:kraft_romano}}

\indent 

We introduce some notation before embarking on the proof.

The density of $P \in \mmp$ wrt a $\sigma$-finite measure $\mu$ is $p$.
The set of densities of all distributions in $\mmp$ is denoted $\boup{p}$.
Similarly, the null and alternative sets of densities are $\boup{p}_0$ and $\boup{p}_1$,
and their union equals $\boup{p}$.
Define $co(\boup{p}')$ to be the convex hull of an arbitrary subset  $\boup{p}' \subseteq \boup{p}$
in a similar fashion as in Equation (\ref{eq:convex}).

The Total Variation (TV) metric between two distributions $P, Q \in \mmp$ with densities $p,q \in \boup{p}$ is defined as
\begin{gather}
d_{TV}(p,q)= \frac{1}{2} \int \left\vert p - q \right\vert ~d\mu.
\end{gather}

\bigskip

The proof of the equivalence of (a) and (b) is shown in three parts.

Part 1: $(a) \Leftrightarrow (a')$
where 
\begin{gather*}
(a): \forall q \in \boup{p}_1 ~\exists \{p_k\}_k \subseteq co(\boup{p}_0) ~\text{such that}~
d_{TV}(p_k,q ) \to 0
\\
(a') : \forall q \in \boup{p}_1 ~ \exists \{p_k\}_k \subseteq co(\boup{p}_0) ~ \text{and} ~ \{ \eps_k \}_k \downarrow 0 
~\text{such that}~ d_{TV}(p_k,q ) < \eps_k ~\forall k
\end{gather*}

Part 1, proof, $(a) \Rightarrow (a')$ :

Fix $q$. For $\eps_k = d_{TV}(p_k,q ) \to 0$, there exists  a monotone sub-sequence
$\eps_{k_j} = d_{TV}(p_{k_j},q ) \downarrow 0$. 
Create new sequences $ \ti{p}_j = p_{k_j}$ and $\ti{\eps}_j = \eps_{k_j}/2 $
so that $d_{TV}(\ti{p}_j ,q )<\ti{\eps}_j $.

Part 1, proof, $(a) \Leftarrow (a')$ : straightforward.

\bigskip

Part 2: $(a') \Leftrightarrow (b')$
where 
\begin{gather*}
(a') : \forall q \in \boup{p}_1 ~ \exists \{p_k\}_k \subseteq co(\boup{p}_0) ~ \text{and} ~ \{ \eps_k \}_k \downarrow 0
~\text{such that}~ d_{TV}(p_k,q ) < \eps_k ~\forall k
\\
(b') : \forall q \in \boup{p}_1 ~ \exists  \{ \eps_k \}_k \downarrow 0
~\text{such that}~ \forall \phi ~  
\int \phi q ~d\mu < \eps_k + \sup_{p \in \boup{p}_0} \int \phi p ~d\mu ~\forall k
\end{gather*}

Part 2, proof, $(a') \Rightarrow (b')$:

Fix $q$, (a') implies there exist sequences $ \{p_k\}_k \subseteq co(\boup{p}_0)$ and $\{ \eps_k \}_k \downarrow 0$
such that $d_{TV}(p_k,q ) < \eps_k$  $\forall k$.
Fix $k$.
Use Theorem \ref{theo:kraft} with $\{ \boup{p}_1 \} = \{ q \}$.
$(a')$ implies

$ \forall \phi ~  
\int \phi q ~d\mu < \eps_k + \sup_{p \in \boup{p}_0} \int \phi p ~d\mu$.

This is true for every $k$ of a sequence $\eps_k$ that converges to zero, given an arbitrary $q$.

Part 2, proof, $(a') \Leftarrow (b')$:

Fix $q$, get $\eps_k$. Fix $k$. 
Use Theorem \ref{theo:kraft} with $\{ \boup{p}_1 \} = \{ q \}$.
$(b')$ implies there exists $p_k \in co(\boup{p}_0)$ such that $d_{TV}(p_k,q)<\eps_k$.
Repeat this for every $k$ to get a sequence $\{p_k\}_k \subseteq co(\boup{p}_0)$
 such that $d_{TV}(p_k,q)<\eps_k ~\forall k$.

\bigskip

Part 3: $(b') \Leftrightarrow (b)$
where 
\begin{gather*}
(b') : \forall q \in \boup{p}_1 ~ \exists  \{ \eps_k \}_k \downarrow 0
~\text{such that}~ \forall \phi ~  
\int \phi q ~d\mu < \eps_k + \sup_{p \in \boup{p}_0} \int \phi p ~d\mu ~\forall k
\\
(b) : \forall \phi ~\text{and}~ q \in \boup{p}_1, 
\int \phi q ~d\mu \leq \sup_{p \in \boup{p}_0 } \int \phi p ~d\mu
\end{gather*}

Part 3, proof $(b') \Rightarrow (b)$:

Fix $q$, get $\eps_k$. Fix $\phi$. It is true that
 
$  
\int \phi q ~d\mu < \eps_k + \sup_{p \in \boup{p}_0} \int \phi p ~d\mu$.

Take limits on both sides,

$  
\int \phi q ~d\mu \leq  \sup_{p \in \boup{p}_0} \int \phi p ~d\mu$.

This is true for every $q$ and every $\phi$.

Part 3, proof $(b') \Leftarrow (b)$:

Straightforward because for arbitrary $\phi$, $q$, and $\{ \eps_k \}_k \downarrow 0$

$  
\int \phi q ~d\mu \leq  \sup_{p \in \boup{p}_0} \int \phi p ~d\mu$

implies that

$  
\int \phi q ~d\mu < \eps_k +  \sup_{p \in \boup{p}_0} \int \phi p ~d\mu$.

%
%
%

$\square$

\subsection{Proof of Theorem \ref{theo:impossible_test}}
\indent

The proof of Theorem \ref{theo:impossible_test} follows the same  lines as the proof of Theorem 1 by \cite{romano2004}
except for the fact that our Assumption \ref{assu:dense} is stated in terms of the LP metric and in terms of the convex hull of $\mmp_0$.

Pick an arbitrary $Q\in \mathbb{P}_{1}$. There exists a
sequence of distributions $\{P_{k}\}_{k=1}^{\infty } \subseteq co(\mmp_0)$ such that $P_k \dto Q$.
Convergence in distribution  is equivalent to $\mathbb{E}_{P_{k}}[g]\rightarrow \mathbb{E}_{Q}[g]$ for every bounded real-valued
function $g$ whose set of discontinuity points has probability zero under $Q$
(Theorem 25.8, \cite{billingsley2008}). In particular, this is true for $g=\phi $ for an arbitrary $\phi $ that is a.s. continuous under $Q$.

Take an arbitrary
sequence $\varepsilon_{n}\rightarrow 0$, and pick a sub-sequence $\{ P_{k_{n}} \}_n$ 
from the sequence $\{P_{k}\}_{k}$ such that
\begin{equation}
-\varepsilon _{n}\leq \mathbb{E}_{Q}\phi -\mathbb{E}_{P_{k_{n}}}\phi \leq
\varepsilon _{n}.
\end{equation}%
Therefore,%
\begin{equation}
\mathbb{E}_{Q}\phi \leq \mathbb{E}_{P_{k_{n}}}\phi +\varepsilon _{n}\leq
\sup_{P\in co(\mmp_{0})}\mathbb{E}_{P}\phi +\varepsilon _{n}.
\end{equation}%
Given $\varepsilon _{n}\rightarrow 0$, it follows that, for $\forall Q\in 
\mathbb{P}$,%
\begin{equation}
\mathbb{E}_{Q}\phi \leq \sup_{P\in co(\mmp_{0})}\mathbb{E}_{P}\phi .
\end{equation}%
Consequently,%
\begin{equation}
\sup_{Q\in \mmp_{1}}\mathbb{E}_{Q}\phi \leq \sup_{P\in co(\mmp_{0})}%
\mathbb{E}_{P}\phi.
\end{equation}%

It is clear that $\sup_{P\in co(\mmp_{0})}\mme_{P}\phi \geq \sup_{P\in \mmp_{0}}\mme_{P}\phi$.
It remains to show that these are equal.
Assume $\sup_{P\in co(\mmp_{0})}\mme_{P}\phi > \sup_{P\in \mmp_{0}}\mme_{P}\phi$.
Select $\eps>0$ small enough such that $\sup_{P\in co(\mmp_{0})}\mme_{P}\phi -\eps > \sup_{P\in \mmp_{0}}\mme_{P}\phi$.
There exists $P_{\eps} \in co(\mmp_{0})$ such that
\begin{equation}
\sup_{P\in co(\mmp_{0})}\mme_{P}\phi \geq \mme_{P_{\eps}} \phi > \sup_{P\in co(\mmp_{0})}\mme_{P}\phi -\eps > \sup_{P\in \mmp_{0}}\mme_{P}\phi.
\end{equation}
By definition, $P_{\eps} = \sum_{i=1}^{N}\alpha_i P_i$ for $N \in \mmn$, $ P_i \in  \mmp_0 ~\forall i$, $\alpha_i \in [0,1] ~\forall i$, and $\sum_{i=1}^N \alpha_i =1$.
Then, $\mme_{P_{\eps}} \phi = \sum_{i=1}^{N}\alpha_i \mme_{P_i} \phi \leq \sup_{P\in \mmp_{0}}\mme_{P}\phi$, a contradiction.
Therefore, $\sup_{P\in co(\mmp_{0})}\mme_{P}\phi$ $ = \sup_{P\in \mmp_{0}}\mme_{P}\phi$ and 
\begin{equation}
\sup_{Q\in \mmp_{1}}\mathbb{E}_{Q}\phi \leq \sup_{P\in \mmp_{0}}\mathbb{E}_{P}\phi.
\end{equation}%

$\square $

\subsection{Proof of Theorem \ref{theo:ci}}

\indent

The proof is a combination of proofs by \cite{dufour1997} and \cite{gleser1987}.

Part (\ref{theo:ci:eq1}):

 Fix $m \in \mu(\mmp)$.
Define $\phi_m = \mmi\{m \not\in C(Z) \}$, and note that 
$\sup\limits_{P\in\mmp(m)}\mme_P \phi_m = \sup\limits_{P\in co(\mmp(m))}\mme_P \phi_m $ (see proof of Theorem \ref{theo:impossible_test}). It follows that  
$1-\alpha \leq \inf\limits_{P\in\mmp(m)} P\left[ m \in C(Z) \right]  
= \inf\limits_{P\in co(\mmp(m))} P\left[ m \in C(Z) \right]$.
Therefore, $\forall P \in co(\mmp(m))$, $P\left[ \mu(P) \in C(Z) \right] \geq 1-\alpha$.

By Assumption \ref{assu:limpt}, there exists $\{P_k\}$ in $co(\mmp(m))$ such that $P_k \dto P^*$.
Then,
\begin{gather}
1-\alpha \leq P_k \left[\mu(P_k) \in C(Z) \right] = P_k \left[m \in C(Z) \right]
\to P^*  \left[m \in C(Z) \right]
\end{gather}
where the convergence follows by  Portmanteau's theorem because $P^*(\partial\{m \in C(Z) \})=0$ (Theorem 29.1 of \cite{billingsley2008}). This proves (\ref{theo:ci:eq1}).

\bigskip

Part (\ref{theo:ci:eq2}):

Pick a sequence $m_n \in \mu(\mmp)$ such that $m_n$ is unbounded. Without loss of generality, assume $m_n \uparrow \infty$. We have that
\begin{gather}
1-\alpha \leq P^*  \left[m_n \in C(Z) \right] \leq P^*  \left[m_n \leq U[C(Z)] \right] .
\end{gather}
Taking the limit  as $n \to \infty$,
\begin{gather}
1-\alpha \leq P^*  \left[ U[C(Z)] = \infty \right] 
\\
\leq P^* \left[U[C(Z)] - L[C(Z)]  = \infty \right]
=  P^* \left[D[C(Z)] = \infty \right].
\end{gather}

\bigskip

Part (\ref{theo:ci:eq3}):

Assumption \ref{assu:limpt} gives a sequence $\{ P_k \}_k$ in $co(\mmp)$
that converges in distribution to $P^*$. 
By assumption, $P^*\left[ \partial \{ D[C(Z)] =\infty ]  \} \right]=0$,
so Portmanteau's theorem gives
$P_k\left[ D[C(Z)] =\infty  ]   \right] \to 
P^* \left[  D[C(Z)] =\infty  ]  \right] \geq 1 - \alpha$.
There exists a sequence $\delta_k \downarrow 0$
such that
$P_k \left[ D[C(Z)] =\infty  ]   \right] \geq 1-\alpha - \delta_k$.

Fix $\eps>0$. 
The set $B_{\eps}(P^*) \cap co(\mmp)$ contains infinitely many $P_k$s from the sequence above. For these  $P_k$s,
\begin{gather}
1- \alpha -\delta_k 
\leq
 P_k \left[ D[C(Z)] =\infty  ] \right] 
 \\
\leq 
\sup_{P \in B_{\eps}(P^*) \cap co(\mmp)} P \left[ D[C(Z)] =\infty  ] \right] 
\\
= \sup_{P \in B_{\eps}(P^*) \cap \mmp} P \left[ D[C(Z)] =\infty  ] \right]
\end{gather}
where the last equality follows by the same argument seen in the proof of (\ref{theo:ci:eq1}) above.
Taking the limit as $k\to \infty$  gives  (\ref{theo:ci:eq3}).

$\square$

\subsection{Coverage of Confidence Intervals - Lemma \ref{lem:covsize}}

\begin{lemma}
\label{lem:covsize} Let $C(Z)$ be constructed as in Equation (\ref{eq:theo:ci_test}). Then,
\begin{equation}
\inf_{P\in \mathbb{P}} {P}\left[ \mu (P)\in C(Z)\right] = 1 - \sup_{m\in \mu (\mathbb{P})} \alpha(m).
\end{equation}%
\end{lemma}

\noindent \textbf{Proof of Lemma \ref{lem:covsize}. }Suppose%
\begin{equation}
\sup_{m\in \mu (\mathbb{P})}\sup_{P\in \mathbb{P}_{0,m}}{P}(\phi
_{m}(Z)=1)=\alpha .
\end{equation}%
Now, pick $\varepsilon >0$. Then, there exists $m_{\varepsilon }$ such that%
\begin{equation}
\alpha -\varepsilon /2\leq \sup_{P\in \mathbb{P}_{0,m_{\varepsilon }}}%
{P}(\phi _{m_{\varepsilon }}(Z)=1)\leq \alpha .
\end{equation}%
There also exists $P_{\varepsilon }\in \mathbb{P}_{0,m_{\varepsilon }}$
such that%
\begin{equation}
\alpha -\varepsilon \leq {P_{\varepsilon }}(\phi _{m_{\varepsilon
}}(Z)=1)\leq \alpha .
\end{equation}%
Rearranging the expression above, we obtain%
\begin{equation}
1-\alpha +\varepsilon \geq {P_{\varepsilon }}(\mu (P_{\varepsilon
})\in C(Z))\geq 1-\alpha .
\end{equation}%
Therefore, we find that%
\begin{equation}
\inf_{P\in \mathbb{P}} {P}\left[ \mu (P)\in C(Z)\right] =1-\alpha ,
\end{equation}%
as we wanted to prove. 

$\square $

\subsection{Proof of Corollary \ref{coro:rdd} }

\indent

Fix $m\in \mathbb{R}$.
Pick an arbitrary $Q\in \mathbb{P}_{1,m}$, and let $m^{\prime }=\mu (Q)\neq
m$. 
Define $g(x) = \mme_Q[Y_i|X_i=x]$.
Construct a sequence of functions $g_{k}:\mathbb{R}\rightarrow \mathbb{R}$, $k=1,2,\ldots $ as follows: 
\begin{equation}
g_{k}(x)=g(x) + (m'-m) \left[ \Lambda\left( {k^2(x-c)} \right) - \mmi\{x \geq c \}\right] 
\end{equation}
where
$\Lambda\left( \cdot \right)$ is the cumulative distribution function (CDF) of the logistic distribution.

The function $g_k$ is infinitely continuously differentiable on $\mmx \setminus \{c\}$, so $g_k \in \m{G} ~\forall k$,
and $\lim\limits_{x \downarrow c}g_k(x) - \lim\limits_{x \uparrow c} g_k(x) = m$.
Moreover, as $k \to \infty$, $g_k(x) \to g(x)$ for every $x \neq c$. Define $P_k$
to be the distribution of $(X_i ~,~Y_i -g(X_i) + g_k(X_i))$ when $(X_i,Y_i) \sim Q$.
It follows that $\mu(P_k)=m $ and $P_k \in \mathbb{P}_{0,m} ~~\forall k$.

It remains to show that $P_k \dto Q$, or equivalently, to show that
\begin{equation}
(X_i ~,~Y_i -g(X_i) + g_k(X_i)) \dto (X_i,Y_i)
\end{equation}
as $k \to \infty$ where $(X_i,Y_i) \sim Q$.
Note that $(X_i ~,~Y_i -g(X_i) + g_k(X_i))= (X_i,Y_i) + (0,g_k(X_i)-g(X_i))$, so it suffices to
show that
$g_k(X_i)-g(X_i) \pto 0$ as  $k \to \infty$.

Define $A_{k}=\{c-k^{-1}< X_{i} < c+k^{-1}\}$, and let $A_{k}^{c}$ be the complement of $A_{k}$.
Fix $\eps>0$.
\begin{gather}
Q\left[ \left|   g_k(X_i)-g(X_i)   \right| > \eps \right]
\\
=Q\left[ \left. \left|   g_k(X_i)-g(X_i)   \right| > \eps ~\right|~ A_k \right] ~ Q\left[ A_k \right] \label{coro:rdd:eq1}
\\
\hspace{4cm}+Q\left[ \left. \left|   g_k(X_i)-g(X_i)   \right| > \eps ~\right|~ A_k^c \right] ~ Q\left[ A_k^c \right]. \label{coro:rdd:eq2}
\end{gather}
Part (\ref{coro:rdd:eq1}) vanishes as $k\to \infty$ by the continuity property of probability measures, 
because $A_{k}\downarrow\{ c \} $ and $Q\left[ \{ c \} \right]=0$ by assumption. 

For part (\ref{coro:rdd:eq2}), note that  $\left|   g_k(x)-g(x)   \right| \leq |m'-m| \Lambda\left( {-k} \right)$ for any $x \in A_k^c$
because 
$\Lambda\left( {k^2(x-c)} \right)$ is strictly increasing in $x$ and symmetric around $x=c$,
so $\left|   g_k(x)-g(x)   \right|$ attains its maximum at $x=c-k^{-1}$ and $x=c+k^{-1}$. Therefore,
\begin{gather}
(\ref{coro:rdd:eq2}) \leq \mmi \left\{ | m' - m | \Lambda\left( {-k} \right) > \eps \right\} ~ Q\left[ A_k^c \right] \to 0
\end{gather}
because $\Lambda\left( {-k} \right) \to 0$ as $k \to \infty$.

Therefore, Assumption \ref{assu:dense} is satisfied for every $m\in \mathbb{R}$. Theorem \ref{theo:impossible_test} applies, and Corollary \ref{coro:ci}
applies with $\mu (\mathbb{P})=\mathbb{R}$. 

$\square $

\subsection{Example of RDD Model with Manipulation}
\label{supp:rdd:manipul}

\indent 

\begin{figure}[H]
\caption{Example of RDD with Manipulation}
\label{fig:rdd:manipul}
\vspace{-.75 \baselineskip}  
\begin{center}
  \begin{minipage}{.5 \textwidth }
    \centering
    (a) Conditional Mean without Manipulation
    
    \includegraphics[width=3.5in]{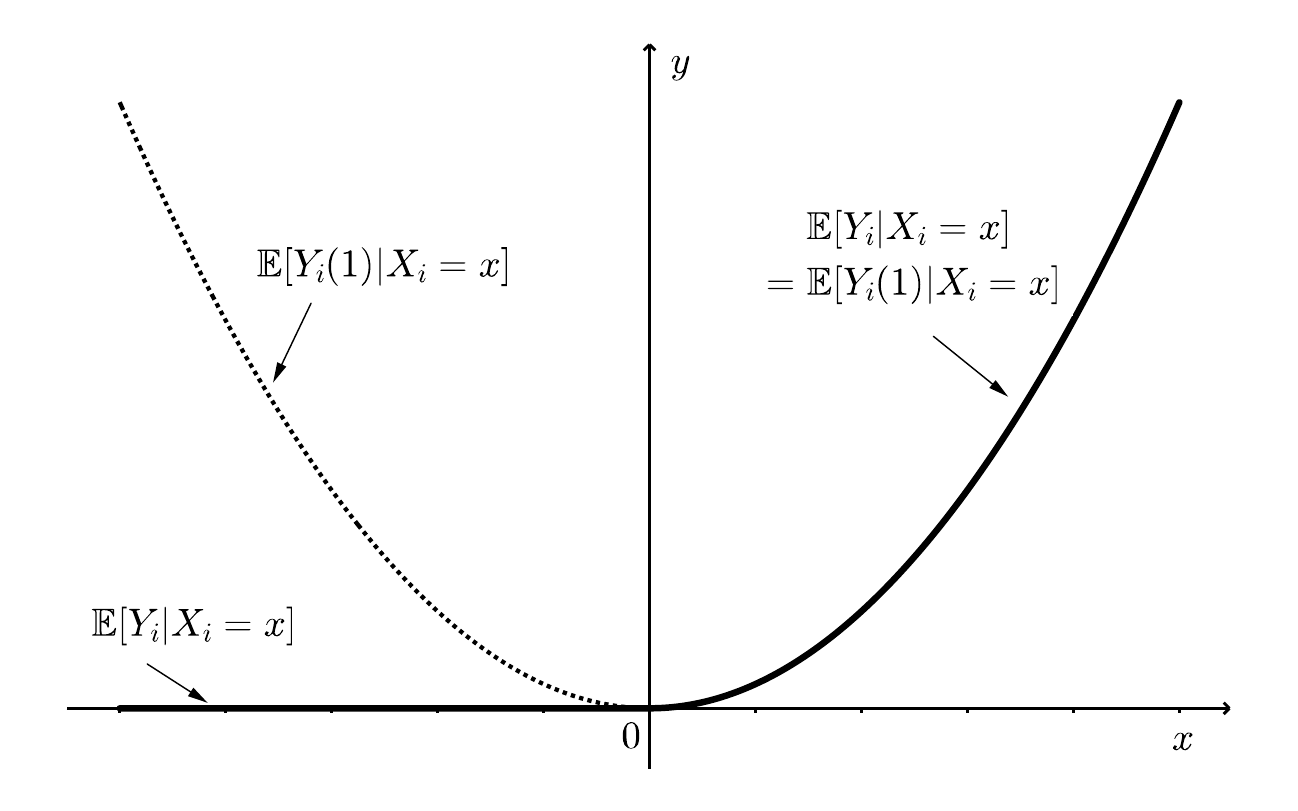}
  \end{minipage}%
  \begin{minipage}{.5 \textwidth }
    \centering
    (b) Conditional Mean with Manipulation

    \includegraphics[width=3.5in]{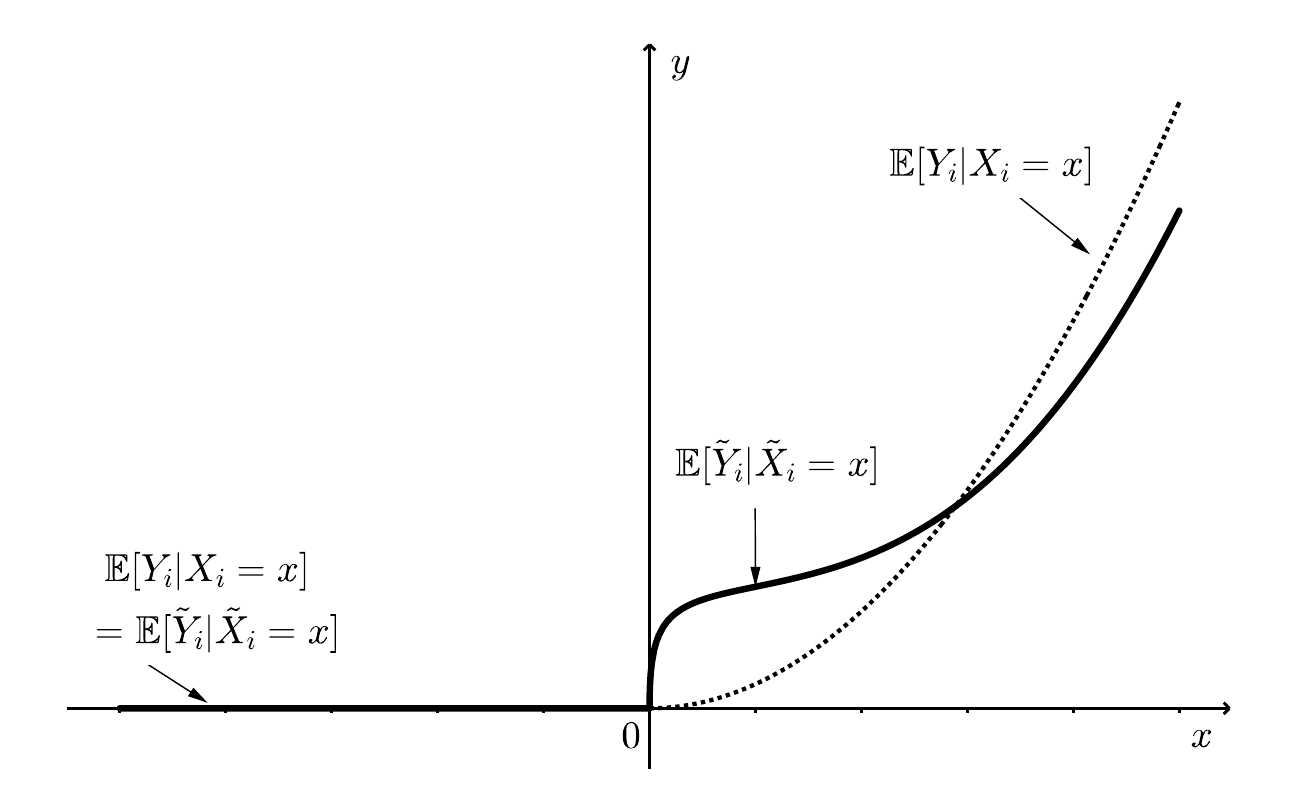}
  \end{minipage}%
  
\end{center}
\vspace{-.5cm}
\caption*{\footnotesize
Notes:
In figure (a) there is no manipulation, and the researcher observes a sample of $(Y_i,X_i)$.
The solid line denotes the conditional mean of the observed outcome given the margin of victory $\mme[Y_i|X_i]$.
The conditional mean of the potential outcome in case of victory
$\mme[Y_i(1)|X_i]$ is the the dotted line, and the potential outcome in case of loss is normalized to zero $Y_i(0)=0$.
In figure (b) there is manipulation, and the researcher observes a sample of $(\tilde Y_i,\tilde X_i)$.
The solid line is $\mme[\tilde Y_i|\tilde X_i]$ while the dotted line depicts $\mme[ Y_i| X_i]$. 
Manipulation increases the slope of the conditional mean function at the cutoff.}
\end{figure}

In this section, we provide an example of DGP with manipulation that gives rise to a model similar to Equation (\ref{dgp:rdd:2}) in our Monte Carlo experiment.
The potential outcome of losing an election is normalized to zero ($Y_i(0)=0$), and the potential outcome of winning an election is
$Y_i(1)$.
Suppose the expected potential gain of winning an election is small for tight elections but large otherwise;
that is, let $\mme[Y_i(1) - Y_i(0) | X_i=x] = \mme[Y_i(1) | X_i=x] = x^2$,
 where $X_i$ is the margin of victory of a given political party in district $i$.
Assume the distribution of $X_i$ is Uniform$[-1,1]$.  Each district is an iid draw $(Y_i(1), X_i, \eps_i)$ from a given distribution,
where $\eps_i$ denotes district $i$'s potential to influence the election outcome in a world where manipulation is possible. 
In a world \textit{without} manipulation, the researcher observes $(Y_i, D_i, X_i)$, where 
$D_i = \mmi \{ X_i \geq 0 \}$ is the victory indicator, and $Y_i = D_i Y_i(1) + (1-D_i)Y_i(0) = D_i Y_i(1)$ is the outcome.
It follows that $\mme[Y_i | X_i = x] = x^2 \mmi\{x \geq 0\}$.
There is no discontinuity at the cutoff, the causal effect is zero, and the DGP is under the null hypothesis of zero effect at the cutoff
(Figure \ref{fig:rdd:manipul}(a)).

\vspace{-1 \baselineskip}

In a world \textit{with} manipulation, the given party in district $i$ decides to influence the election if the expected potential gain of doing so is positive.
In other words, if the margin of victory without manipulation $X_i$ leads to the loss of the election, 
and the expected potential gain of winning the election $\mme[Y_i(1)-Y_i(0)|X_i]$ is strictly positive, then 
the party decides to manipulate.
Thus, manipulation occurs if $X_i<0$, and the margin of victory changes from $X_i$ to $\eps_i>0$.
Although the party manipulates as little as possible to win the election, it does not have perfect control over its vote share. 
Assume $\eps_i=5 \chi^2_{(3)}$, that is, five times a Chi-square distribution with three degrees of freedom.
The pdf $f_{\eps}$ evaluated at zero equals zero, but it is highly sloped to the right of zero.
Let $\tilde{X}_i$ be the manipulated margin of victory defined as 
$\tilde{X}_i = \mmi\{X_i < 0\} \eps_i + \mmi\{X_i \geq 0  \} X_i$. 
The researcher observes $(\tilde{Y}_i, \tilde{D}_i, \tilde{X}_i)$, where $\tilde{Y}_i = \tilde{D}_i Y_i(1) + (1-\tilde{D}_i)Y_i(0) = \tilde{D}_i Y_i(1)$,
and $\tilde{D}_i = \mmi \{ \tilde{X}_i \geq 0 \}$ is the victory indicator.
The conditional mean function under manipulation is given by 
\begin{gather*}
\mme[\tilde{Y}_i | \tilde{X}_i = x] = \mmi\{x \geq 0\} \mme[Y_i(1) | \tilde{X}_i=x]
\\
=\mmi\{x \geq 0\} \mme\Big[Y_i(1)  ~ \Big|~ \{ \eps_i=x, X_i < 0 \} \text{ or } \{ X_i=x, X_i \geq  0\} \Big]
\\
=\mmi\{x \geq 0\} \Big\{ \theta(x)  \mme\left[Y_i(1)  ~ |~ \eps_i=x, X_i < 0 \right] + (1-\theta(x)) \mme\left[Y_i(1)  ~ |~ X_i=x, X_i \geq  0 \right] \Big\}
\\
=\mmi\{x \geq 0\} \Big\{ \theta(x)  \mme\left[Y_i(1)  ~ |~ X_i < 0 \right] + (1-\theta(x)) \mme\left[Y_i(1)  ~ |~ X_i=x \right] \Big\}
\\
=\mmi\{x \geq 0\} \Big\{ \theta(x)  (1/3) +
(1-\theta(x)) x^2 \Big\},
\end{gather*}
where the weight
\begin{gather*}
\theta(x) = \frac{f_{\eps}(x) \mmp(X_i < 0)}{f_{\eps}(x) \mmp(X_i < 0) + f_X(x) \mmi\{x \geq 0\} }
\\
= \frac{f_{\eps}(x) 0.5 }{f_{\eps}(x) 0.5 + 0.5 \mmi\{x \geq 0\} }
\end{gather*}
is such that $\theta(x)\in[0,1]$, $\theta(0)=0$, $\theta(x)$ is continuous in $x$, and it is positively and highly sloped near $x=0$.  
The conditional mean function
$\mme[\tilde{Y}_i | \tilde{X}_i = x]$ increases sharply at the cutoff (Figure \ref{fig:rdd:manipul}(b))
because districts with low $X_i$ and high causal effects manipulate their $X_i$ to $\tilde{X}_i=\eps_i$ to the right of the cutoff. 
There is no discontinuity at the cutoff, and the DGP is still under the null hypothesis of zero effect.
However, manipulation makes it harder to distinguish a zero effect from a positive effect at the cutoff.

\subsection{Proof of Corollary \ref{coro:bun}}

\indent

 Fix $Q \in \mmp_1$ with CDF $F_Q(x)$. The CDF $F_Q(x)$ has a jump discontinuity of size $\delta>0$ at $x=x_0$.
 Call $f_Q$ the derivative of $F_Q$ at $x\neq x_0$, which is a continuous function of $x$ for every $x\neq x_0$.
 The integral of $f_Q$ over $\mmr$ equals $1-\delta$.
 The side limits of $f_Q$ at $x_0$, $f_Q(x_0^+)$ and $f_Q(x_0^-)$, may be different from each other.
Pick a sequence $\eps_k \downarrow 0$. 
Construct a continuous ``hat-shaped'' function $g_k(x):[x_0-\eps_k;x_0+\eps_k] \to \mmr$ such that: (i) $g_k(x_0-\eps_k)=f_Q(x_0-\eps_k)$;
 (ii) $g_k(x_0 + \eps_k)=f_Q(x_0 + \eps_k)$; (iii) $g_k(x)$ has constant and positive slope for $x \leq x_0$,
 and constant and negative slope for $x \geq x_0 $; (iv) $g_k(x) \geq f_Q(x)$; and (v) $\int (g_k(x) - f_Q(x)) ~dx = \delta $.
It is always possible to construct such a function for a small enough $\eps_k$.
Define $f_{P_k}(x) = f_Q(x) + \mmi\{x_0-\eps_k \leq x \leq  x_0+\eps_k \} \left( g_k(x) - f_Q(x) \right)$.
This is a continuous PDF function, and let it define the distribution $P_k$.
Then the CDF $F_{P_k}$ converges to $F_Q$ as $k \to \infty$ at every continuity point of $F_Q$,
so that $P_k \dto Q $. 

$\square$

\subsection{Proof of Corollary \ref{coro:exo}}

\indent

 Fix $m\in \mathbb{R}$.
Choose an arbitrary $Q\in \mathbb{P}_{1,m}$, and let $m^{\prime }=\mu (Q)\neq m$.
Define $g(x,w) = \mme_Q[Y_i|X_i=x,W_i=w]$,
and $\tau_Q(w)= g(x,w) - \lim_{x \downarrow 0}g(x,w)$.

Construct a sequence of functions $g_{k}:\mmx \times \mmw \rightarrow \mathbb{R}$, $k=1,2,\ldots $ as follows: 
\begin{equation}
g_{k}(x,w)=g(x,w) + (\tau_Q(w) -m) \left[ \mmi\{x > 0 \} - \Lambda\left( k^2 x \right) \right] 
\end{equation}
where
$\Lambda\left( \cdot \right)$ is the CDF of the logistic distribution.

The function $g_k$ is infinitely many times continuously differentiable wrt $x$  on $\{ \mmx \setminus \{c\}\} \times \mmw$, so $g_k \in \m{G} ~\forall k$.
Also, $g_k(0,w) - \lim\limits_{x \downarrow 0}g_k(x,w) = m$.
Moreover, as $k \to \infty$, $g_k(x,w) \to g(x,w)$ pointwise. Define $P_k$
to be the distribution of $(X_i ~,~ W_i ~,~Y_i -g(X_i,W_i) + g_k(X_i,W_i))$ when $(X_i,W_i,Y_i) \sim Q$.
It follows that $\mu(P_k)=m $ and $P_k \in \mathbb{P}_{0,m} ~~\forall k$.

It remains to show that $P_k \dto Q$, or equivalently, to show that
\begin{equation}
(X_i ~,~ W_i ~,~Y_i -g(X_i,W_i) + g_k(X_i,W_i)) \dto (X_i,W_i,Y_i)
\end{equation}
as $k \to \infty$ where $(X_i,Y_i) \sim Q$.
Note that $(X_i ~,~ W_i ~,~Y_i -g(X_i,W_i) + g_k(X_i,W_i))= (X_i,W_i,Y_i) + (0,g_k(X_i,W_i)-g(X_i,W_i))$, so it suffices to
show that
$g_k(X_i,W_i)-g(X_i,W_i) \pto 0$ as  $k \to \infty$.

Define $A_{k}=\{0 < X_{i} < k^{-1} \}$, and let $A_{k}^{c}$ be the complement of $A_{k}$.
Fix $\eps>0$.
\begin{gather}
Q\left[ \left|   g_k(X_i,W_i)-g(X_i,W_i)   \right| > \eps \right]
\\
=Q\left[ \left. \left|   g_k(X_i,W_i)-g(X_i,W_i)   \right| > \eps ~\right|~ A_k \right] ~ Q\left[ A_k \right] \label{coro:exo:eq1}
\\
\hspace{4cm}+Q\left[ \left. \left|   g_k(X_i,W_i)-g(X_i,W_i)   \right| > \eps ~\right|~ A_k^c \right] ~ Q\left[ A_k^c \right]. \label{coro:exo:eq2}
\end{gather}
Part (\ref{coro:exo:eq1}) vanishes as $k\to \infty$ by the continuity property of probability measures 
because $A_{k}\downarrow\{ \emptyset \} $ where $\emptyset$ denotes the empty set and has zero probability. 

For part (\ref{coro:exo:eq2}),   $\left|   g_k(x,w)-g(x,w)   \right| \leq | \tau_Q(w) -m | | 1 - \Lambda\left( {k} \right) |$ for any $w$ and any $x \in A_k^c$
because 
$1- \Lambda\left( k^2 x \right)$ is strictly decreasing in $x$.
For fixed $w$, $\left|   g_k(x,w)-g(x,w)   \right|$ attains its maximum at $x=k^{-1}$. Therefore,
\begin{gather}
(\ref{coro:exo:eq2}) \leq \mmp \left\{  | \tau_Q(W_i) -m | | 1 - \Lambda\left( {k} \right) | > \eps \right\} ~ Q\left[ A_k^c \right] \to 0
\end{gather}
because $\Lambda\left(  k \right) \to 1$ as $k \to \infty$ and $| \tau_Q(W_i) -m |$ is bounded.

$\square $

\subsection{Simulations - RDD}
\label{supp:simul:rdd}

\indent 

This section contains additional results of the RDD simulation in the main text. 
The size and power analyses in the main text use the 5\% nominal level. 
This section has the same analyses using 1\% and 10\% nominal levels.
It also has the simulated critical values under the various choices of null $(\tau,M)$-models.

\begin{landscape}
  
\begin{table}
\begin{center}
\caption{Simulated Rejection Rates and Critical Values}
\label{supp:table:rdd:simul}
\scriptsize

{\normalsize {Panel 1:} \textit{Rejection Rate under the Null Model $(\tau,M)$ Using Critical Values Simulated under Model $(\tau,0)$ }}

	\begin{minipage}{9cm}
    \centering
    (a) Nominal Size 1\% 
    
     \begin{tabular}{c  c  c  c  c  c  c  } 
 \hline \hline 
$\tau$  & $ M= 0 $  & $ M= 2 $  & $ M= 4 $  & $ M= 6 $  & $ M= 8 $  & $ M= 10 $ 
\\
\hline
 .01  &   0.0100 &   0.0097 &   0.0108 &   0.0142 &   0.0108 &   0.0123 \\ 
 .02  &   0.0100 &   0.0130 &   0.0147 &   0.0120 &   0.0139 &   0.0132 \\ 
 .03  &   0.0100 &   0.0195 &   0.0192 &   0.0212 &   0.0210 &   0.0224 \\ 
 .04  &   0.0100 &   0.0231 &   0.0257 &   0.0232 &   0.0243 &   0.0276 \\ 
 .05  &   0.0100 &   0.0274 &   0.0323 &   0.0335 &   0.0373 &   0.0345 \\ 
 .06  &   0.0100 &   0.0339 &   0.0467 &   0.0523 &   0.0591 &   0.0577 \\ 
 .07  &   0.0100 &   0.0402 &   0.0546 &   0.0612 &   0.0702 &   0.0737 \\ 
 .08  &   0.0100 &   0.0451 &   0.0674 &   0.0827 &   0.0917 &   0.0968 \\ 
\hline  \hline 
\end{tabular} 

  	\end{minipage} 
  	\hspace{.5cm}
	\begin{minipage}{9cm }
    \centering
    (b) Nominal Size 10\% 
    
     \begin{tabular}{c  c  c  c  c  c  c  } 
 \hline \hline 
$\tau$  & $ M= 0 $  & $ M= 2 $  & $ M= 4 $  & $ M= 6 $  & $ M= 8 $  & $ M= 10 $ 
\\
\hline
 .01  &   0.1000 &   0.1096 &   0.1111 &   0.1170 &   0.1159 &   0.1127 \\ 
 .02  &   0.1000 &   0.1253 &   0.1309 &   0.1270 &   0.1337 &   0.1254 \\ 
 .03  &   0.1000 &   0.1587 &   0.1618 &   0.1651 &   0.1629 &   0.1707 \\ 
 .04  &   0.1000 &   0.1733 &   0.1803 &   0.1818 &   0.1827 &   0.1901 \\ 
 .05  &   0.1000 &   0.1801 &   0.2049 &   0.2097 &   0.2244 &   0.2177 \\ 
 .06  &   0.1000 &   0.2126 &   0.2561 &   0.2635 &   0.2810 &   0.2805 \\ 
 .07  &   0.1000 &   0.2172 &   0.2596 &   0.2867 &   0.3022 &   0.3072 \\ 
 .08  &   0.1000 &   0.2244 &   0.2957 &   0.3238 &   0.3427 &   0.3594 \\ 
\hline  \hline 
\end{tabular} 

  	\end{minipage} 


{\normalsize {Panel 2:} \textit{Rejection Rate under the Alternative Model $(\tau,\infty)$ Using Critical Values Simulated under Model $(\tau,M)$ }}

	\begin{minipage}{9cm}
    \centering
    (a) Nominal Size 1\% 
    
     \begin{tabular}{c  c  c  c  c  c  c  } 
 \hline \hline 
$\tau$  & $ M= 0 $  & $ M= 2 $  & $ M= 4 $  & $ M= 6 $  & $ M= 8 $  & $ M= 10 $ 
\\
\hline
 .01  &   0.0113 &   0.0104 &   0.0099 &   0.0100 &   0.0100 &   0.0101 \\ 
 .02  &   0.0159 &   0.0111 &   0.0105 &   0.0102 &   0.0101 &   0.0103 \\ 
 .03  &   0.0223 &   0.0111 &   0.0112 &   0.0109 &   0.0111 &   0.0105 \\ 
 .04  &   0.0286 &   0.0133 &   0.0114 &   0.0109 &   0.0114 &   0.0108 \\ 
 .05  &   0.0390 &   0.0188 &   0.0137 &   0.0125 &   0.0114 &   0.0111 \\ 
 .06  &   0.0668 &   0.0211 &   0.0149 &   0.0132 &   0.0124 &   0.0114 \\ 
 .07  &   0.0869 &   0.0237 &   0.0150 &   0.0139 &   0.0128 &   0.0132 \\ 
 .08  &   0.1165 &   0.0295 &   0.0184 &   0.0164 &   0.0132 &   0.0130 \\ 
\hline  \hline 
\end{tabular} 

  \end{minipage} 
  \hspace{.5cm}
    \begin{minipage}{9cm}
    \centering
    (b) Nominal Size 10\% 
    
     \begin{tabular}{c  c  c  c  c  c  c  } 
 \hline \hline 
$\tau$  & $ M= 0 $  & $ M= 2 $  & $ M= 4 $  & $ M= 6 $  & $ M= 8 $  & $ M= 10 $ 
\\
\hline
 .01  &   0.1150 &   0.1013 &   0.1004 &   0.1003 &   0.1001 &   0.1001 \\ 
 .02  &   0.1411 &   0.1058 &   0.1026 &   0.1025 &   0.1013 &   0.1013 \\ 
 .03  &   0.1695 &   0.1130 &   0.1060 &   0.1045 &   0.1031 &   0.1030 \\ 
 .04  &   0.2003 &   0.1239 &   0.1129 &   0.1078 &   0.1064 &   0.1038 \\ 
 .05  &   0.2389 &   0.1361 &   0.1182 &   0.1128 &   0.1116 &   0.1079 \\ 
 .06  &   0.3093 &   0.1543 &   0.1319 &   0.1203 &   0.1154 &   0.1109 \\ 
 .07  &   0.3335 &   0.1792 &   0.1438 &   0.1286 &   0.1192 &   0.1161 \\ 
 .08  &   0.4073 &   0.2085 &   0.1567 &   0.1342 &   0.1269 &   0.1200 \\ 
\hline  \hline 
\end{tabular} 

  \end{minipage} 


{\normalsize {Panel 3:} \textit{Critical Values Simulated under Null Model $(\tau,M)$} }

  \begin{minipage}{9cm}
    \centering
    (a) Nominal Size 1\% 
    
     \begin{tabular}{c  c  c  c  c  c  c  } 
 \hline \hline 
$\tau$  & $ M= 0 $  & $ M= 2 $  & $ M= 4 $  & $ M= 6 $  & $ M= 8 $  & $ M= 10 $ 
\\
\hline
 .01  &   3.0222 &   3.0074 &   3.0631 &   3.1741 &   3.0641 &   3.0927 \\ 
 .02  &   3.0669 &   3.1787 &   3.2466 &   3.1459 &   3.1743 &   3.2177 \\ 
 .03  &   3.0506 &   3.3975 &   3.4213 &   3.4238 &   3.3659 &   3.4641 \\ 
 .04  &   3.0810 &   3.5389 &   3.5381 &   3.4846 &   3.4706 &   3.5916 \\ 
 .05  &   3.0786 &   3.5443 &   3.6329 &   3.6819 &   3.7531 &   3.7461 \\ 
 .06  &   2.9715 &   3.5573 &   3.7476 &   3.8374 &   3.9043 &   3.8987 \\ 
 .07  &   2.9829 &   3.7526 &   3.9103 &   3.9917 &   3.9896 &   4.0030 \\ 
 .08  &   2.9536 &   3.7593 &   4.0283 &   4.1033 &   4.2614 &   4.2798 \\ 
\hline  \hline 
\end{tabular} 

  \end{minipage}%
  \hspace{.5cm}
 \begin{minipage}{9cm}
    \centering
    (b) Nominal Size 10\% 
    
     \begin{tabular}{c  c  c  c  c  c  c  } 
 \hline \hline 
$\tau$  & $ M= 0 $  & $ M= 2 $  & $ M= 4 $  & $ M= 6 $  & $ M= 8 $  & $ M= 10 $ 
\\
\hline
 .01  &   1.8700 &   1.9312 &   1.9323 &   1.9577 &   1.9602 &   1.9437 \\ 
 .02  &   1.8801 &   2.0155 &   2.0295 &   2.0221 &   2.0558 &   2.0240 \\ 
 .03  &   1.8466 &   2.1393 &   2.1621 &   2.1594 &   2.1638 &   2.2066 \\ 
 .04  &   1.8822 &   2.2410 &   2.2767 &   2.2879 &   2.2917 &   2.3270 \\ 
 .05  &   1.9056 &   2.3180 &   2.4002 &   2.4403 &   2.4903 &   2.4722 \\ 
 .06  &   1.8316 &   2.3812 &   2.5359 &   2.5594 &   2.6196 &   2.6234 \\ 
 .07  &   1.8772 &   2.4403 &   2.5859 &   2.7054 &   2.7526 &   2.7885 \\ 
 .08  &   1.8651 &   2.4662 &   2.7184 &   2.8357 &   2.8925 &   2.9366 \\ 
\hline  \hline 
\end{tabular} 

  \end{minipage}%

  \begin{minipage}{9cm}
    \centering
    (c) Nominal Size  5\% 

      \begin{tabular}{c  c  c  c  c  c  c  } 
 \hline \hline 
$\tau$  & $ M= 0 $  & $ M= 2 $  & $ M= 4 $  & $ M= 6 $  & $ M= 8 $  & $ M= 10 $ 
\\
\hline
 .01  &   2.2543 &   2.2888 &   2.3158 &   2.3511 &   2.3279 &   2.3525 \\ 
 .02  &   2.2491 &   2.4041 &   2.4038 &   2.3951 &   2.4190 &   2.4311 \\ 
 .03  &   2.1983 &   2.5607 &   2.5687 &   2.5776 &   2.5385 &   2.6150 \\ 
 .04  &   2.2350 &   2.6669 &   2.6713 &   2.6819 &   2.7055 &   2.7515 \\ 
 .05  &   2.2390 &   2.7010 &   2.8138 &   2.8575 &   2.9340 &   2.8821 \\ 
 .06  &   2.2030 &   2.7688 &   2.9399 &   2.9920 &   3.0553 &   3.0453 \\ 
 .07  &   2.2713 &   2.8591 &   3.0370 &   3.0936 &   3.1734 &   3.2180 \\ 
 .08  &   2.2556 &   2.9023 &   3.1370 &   3.2685 &   3.3442 &   3.3817 \\ 
\hline  \hline 
\end{tabular} 

  \end{minipage}%

\caption*{\footnotesize
Notes: Model $(\tau,M)$ refers to Equation (\ref{dgp:rdd:2}) in the main text. 
The estimates for the Wald test are obtained by the robust bias-corrected method of 
\cite{cattaneo2014calonico} and implemented using the STATA package ``rdrobust''.
}

\end{center}

\end{table}

\end{landscape}

\end{singlespace}

\end{document}